\documentclass{SciPost}

\usepackage{subcaption}
\usepackage{graphicx} 
\usepackage[utf8]{inputenc} 
\usepackage[T1]{fontenc}    
\usepackage{hyperref}       
\usepackage{url}            
\usepackage{booktabs}       
\usepackage{amsmath}
\usepackage{nicefrac}       
\usepackage{microtype}      
\usepackage{doi}
\usepackage{bm}
\usepackage{bbold}
\usepackage{comment}
\usepackage{wasysym}
\usepackage[dvipsnames]{xcolor} 
\usepackage{amsthm}
\usepackage{parskip}
\usepackage{tcolorbox}

\usepackage{algorithm}
\usepackage{algorithmic}

\newcommand{\dd}{\mathrm{d}}

\newtcolorbox{propbox}{
  colback=Emerald!5!white,
  colframe=Emerald!40!black,
  boxrule=0.5pt,
  arc=3pt,
  left=3pt,
  right=3pt,
  top=0pt,
  bottom=0pt
}

\newcommand{\eps}{\varepsilon}

\newtheorem{remark}{Remark}[section]

\hypersetup{
    colorlinks,
    linkcolor={red!50!black},
    citecolor={blue!50!black},
    urlcolor={blue!80!black}
}

\usepackage[bitstream-charter]{mathdesign}
\DeclareSymbolFont{usualmathcal}{OMS}{cmsy}{m}{n}
\DeclareSymbolFontAlphabet{\mathcal}{usualmathcal}

\fancypagestyle{SPstyle}{
\fancyhf{}
\lhead{\colorbox{scipostblue}{\bf \color{white} ~SciPost Physics Lecture Notes }}
\rhead{{\bf \color{scipostdeepblue} ~Submission }}

\fancyfoot[C]{\textbf{\thepage}}
}

\begin{document}

\pagestyle{SPstyle}

\begin{center}{\Large \textbf{\color{scipostdeepblue}{
Large Deviations in Stochastic Dynamics:\\ A Hamiltonian Perspective
}}}\end{center}

\begin{center}\textbf{
Lorenzo Cirigliano\textsuperscript{1,2,$\star$},
William D. Pi\~{n}eros\textsuperscript{3,$\dagger$}, 
Anna Poggialini\textsuperscript{2,$\ddagger$}, \\
Grant M. Rotskoff\textsuperscript{4,$\spadesuit$} and 
Eric Vanden-Eijnden\textsuperscript{5,6$\clubsuit$,$\bullet$}
}\end{center}

\begin{center}
{\bf 1} The Abdus Salam International Centre for Theoretical Physics ICTP, Trieste, Italy
\\
{\bf 2} Department of Physics, La Sapienza University of Rome, Italy
\\
{\bf 3} Department of Physics and Material Science, University of Luxembourg, Luxembourg
\\
{\bf 4} Department of Chemistry, Stanford University, USA
\\
{\bf 5} Courant Institute of Mathematical Sciences, New York University, USA
\\
{\bf 6} Machine Learning Lab, Capital Fund Management, Paris, France
\\[\baselineskip]
$\star$ \href{mailto:email1}{\small lorenzo.cirigliano@uniroma1.it}\,,\quad
$\dagger$ \href{mailto:email2}{\small william.pineros@unicaen.fr}\,,\\
$\ddagger$ \href{mailto:email3}{\small anna.poggialini@uniroma1.it}\,,\quad
$\spadesuit$ \href{mailto:email4}{\small rotskoff@stanford.edu}\,,\\
$\clubsuit$ \href{mailto:email4}{\small eve2@nyu.edu, eric.vanden-eijnden@cfm.com},\quad\\
$\bullet${Authors are listed in alphabetical order.}
\end{center}

\section*{\color{scipostdeepblue}{Abstract}}
\textbf{%
These lecture notes present an introduction to large deviations in stochastic dynamical systems, adapted from a lecture series given at the Les Houches school ``Large Deviations: Theory and Applications'' in July 2024. The exposition develops the theory through a unified methodology: WKB analysis of backward Kolmogorov equations, which reveals the Hamiltonian structure underlying all large deviation principles considered here. %
We treat three asymptotic regimes systematically. First, the small-noise limit for stochastic differential equations and the large-volume limit for Markov jump processes, where the Large Deviation Principle emerges from the Law of Large Numbers through progressively finer characterizations of fluctuations. Second, the long-time limit combined with small noise, yielding the quasipotential formalism for metastability and transition rates in both equilibrium and nonequlibrium systems. Third, the long-time limit at finite noise, connecting to spectral methods via tilted generators and enabling reduction of multiscale systems. %
We survey computational approaches including minimum action methods, their geometric variants for quasipotential calculations, and neural network techniques. The Schl\"ogl model of chemical bistability serves as a running example, illustrating both scalar and spatially extended settings. Appendices provide self-contained background on stochastic calculus and Markov jump processes.
}

\vspace{\baselineskip}

\noindent\textcolor{white!90!black}{%
\fbox{\parbox{0.975\linewidth}{%
\textcolor{white!40!black}{\begin{tabular}{lr}%
  \begin{minipage}{0.6\textwidth}%
    {\small Copyright attribution to authors. \newline
    This work is a submission to SciPost Physics Lecture Notes. \newline
    License information to appear upon publication. \newline
    Publication information to appear upon publication.}
  \end{minipage} & \begin{minipage}{0.4\textwidth}
    {\small Received Date \newline Accepted Date \newline Published Date}%
  \end{minipage}
\end{tabular}}
}}
}


\vspace{10pt}
\noindent\rule{\textwidth}{1pt}
\tableofcontents
\noindent\rule{\textwidth}{1pt}
\vspace{10pt}

\section{Introduction}

Rare events---nucleation of a new phase, switching between metastable states, extreme fluctuations of currents or densities---often have dramatic consequences for the behavior of stochastic systems. Yet because they occur with probabilities that are exponentially small on the natural time scales of the dynamics, they are inaccessible to direct simulation or standard perturbation theory. Large deviation theory provides the natural mathematical framework for analyzing such events: it characterizes exponentially rare fluctuations through variational principles that transform probabilistic questions into deterministic optimization problems. These notes, adapted from lectures given at the Les Houches school ``Large Deviations: Theory and Applications'' (July 2024), develop this theory for stochastic dynamical systems with an emphasis on the computational structures that emerge.

The unifying thread of these notes is methodological: in each setting we derive the Large Deviation Principle (LDP) from a multiscale expansion technique, namely the WKB (Wentzel--Kramers--Brillouin) ansatz. We apply this formal asymptotic expansion to the appropriate backward equation to obtain the LDP. This approach is simple and transparent, revealing the Hamiltonian structure that underlies all large deviation problems in a unified way.
This Hamiltonian structure reflects the G\"artner--Ellis theorem at the infinitesimal level: the Hamiltonian is a cumulant generating function, and the Lagrangian---obtained via Legendre transform---is the corresponding rate function. The WKB approach is also generic, applying with minimal modification across limits taken in noise, spatial discretization, and time. While these formal derivations are not rigorous, the resulting statements can be established rigorously using the theory of viscosity solutions or probabilistic methods; we refer to the foundational texts~\cite{FreidlinWentzell2012,dembo2009large,varadhan1984large,denhollander2000,fengkurtz2006} for complete proofs.

The notes are organized around two fundamental asymptotic regimes---small noise and long times---and their interplay.

The first two sections develop the LDP at fixed, finite time. Section~\ref{sec:LDP_trajectories} treats stochastic differential equations (SDEs) in the small-noise limit ($\eps \to 0$). We present the succession from the Law of Large Numbers to the Central Limit Theorem to the Large Deviation Principle, showing how each provides a progressively richer characterization of the stochastic dynamics: typical behavior, Gaussian fluctuations, and rare large deviations, respectively. The LDP transforms the evaluation of rare-event probabilities into a variational problem with natural Hamiltonian structure. Section~\ref{Sec.MJP} extends this framework to Markov jump processes (MJPs) in the large-volume limit ($N \to \infty$), where the small parameter $\eps = 1/N$ arises as a continuum limit for spatial discretization. We illustrate the theory through the Schl\"ogl model of chemical bistability and discrete diffusion on a ring.

The next two sections address the long-time limit ($T \to \infty$). Section~\ref{Sec.QP} develops the quasipotential formalism, which characterizes metastability and transition rates between stable states. The quasipotential emerges naturally in the combined limit of small noise (or large volume) and long time, and provides a principled generalization of potential landscapes to nonequlibrium systems. We illustrate the formalism through both single-compartment and spatially extended versions of the Schl\"ogl model. Section~\ref{sec:LDPtime} treats the long-time limit with $\mathcal{O}(1)$ noise, where the LDP characterizes fluctuations of time-averaged observables. We show how this regime connects to spectral methods via tilted generators, and how it enables systematic reduction of multiscale systems through averaging over fast degrees of freedom.

Finally, Section~\ref{sec:numerical} surveys numerical methods for computing large deviation quantities across all these settings. We discuss both classical approaches based on minimum action methods and modern machine-learning techniques, illustrating the latter through applications to entropy production in active matter systems.

\paragraph{Reading guide:} Readers seeking a first pass through the material may focus on Sections~\ref{sec:LDP_trajectories} and~\ref{Sec.QP}, which develop the core theory for SDEs and the quasipotential, returning to Markov jump processes (Section~\ref{Sec.MJP}) and tilted generators (Section~\ref{sec:LDPtime}) afterward. Section~\ref{sec:numerical} on computational methods can be read independently once Section~\ref{Sec.QP} is understood. Appendices~\ref{appendix:ITO} and \ref{appendix:MJP} provide background on stochastic differential equations and Markov jump processes; readers new to the subject should consult them before reading the main text.

\section{Small-noise limit: Stochastic differential equations}
\label{sec:LDP_trajectories}

In this section, we present the theory of large deviations for stochastic differential equations (SDEs) in the limit of small noise. We establish the basic definitions, introduce the backward Kolmogorov equation, and show how the Large Deviation Principle emerges naturally from a WKB analysis, leading to a Hamiltonian formulation that extends to general diffusion processes.

Consider a general stochastic differential equation (SDE)
\begin{equation}
    \dd X_t^{\eps} = b(t, X_t^{\eps})\,\dd t + \sqrt{2\eps}\, \sigma(t, X_t^{\eps})\, \dd W_t,
    \label{eq:SDE_definition}
\end{equation}
with initial condition $X_{t=0}^{\eps}=x_0$. Here $X_t^{\eps} \in \mathbb{R}^d$ denotes the state variables, $b: \mathbb{R}\times\mathbb{R}^d \to \mathbb{R}^d$ is the deterministic drift, $\sigma: \mathbb{R}\times\mathbb{R}^d \to \mathbb{R}^{d\times d}$ is the volatility, and $W_t$ is a $d$-dimensional Wiener process. The parameter $\eps > 0$ controls the noise strength, and we are interested in the limit $\eps \downarrow 0$. The noise covariance matrix is $D(t,x) = (\sigma \sigma^{\top})(t,x)$, which we assume positive-definite (and hence invertible). For a review of basic properties of SDEs we refer the reader to Appendix~\ref{appendix:ITO}.

\subsection{The small-noise limit: three levels of description}

The limit $\eps \to 0$ admits three progressively more information-rich levels of description, corresponding to the Law of Large Numbers, Central Limit Theorem, and Large Deviation Principle.

\subsubsection{Law of Large Numbers (LLN)}

In the absence of noise ($\eps=0$), the dynamics becomes deterministic. Let $X^{0}_t$ denote the solution of the ordinary differential equation
\begin{equation}
    \dot{X}_t = b(t, X_t), \qquad X_0 = x_0.
\end{equation}
The Law of Large Numbers states that
\begin{equation}
    \lim_{\eps \downarrow 0} \mathbb{E}^{x_0}\left[ \sup_{0 \leq t \leq T} \left|X^{\eps}_{t} - X_{t} \right| \right] = 0, \quad \forall\, T>0,
\end{equation}
where $\mathbb{E}^{x_0}[\cdot]$ denotes the conditional expectation over realizations of the stochastic process starting from $X^\eps_{t=0}=x_0$. Equivalently, for any $\delta > 0$,
\begin{equation}
    \lim_{\eps \downarrow 0} \mathbb{P}^{x_0}\left[\sup_{0 \leq t \leq T} \left|X^{\eps}_{t} - X_{t} \right| > \delta \right] = 0.
\end{equation}
Thus, typical trajectories of the SDE converge uniformly to the deterministic trajectory as $\eps \downarrow 0$.

\subsubsection{Central Limit Theorem (CLT)}

The CLT characterizes the distribution of fluctuations around the deterministic trajectory on the natural scale $\sqrt{\eps}$. Define the rescaled deviation
\begin{equation}
    \eta_t^{\eps} = \frac{X_t^{\eps} - X_t}{\sqrt{\eps}}.
\end{equation}
Differentiating and expanding $b(t, X_t + \sqrt{\eps}\eta_t^{\eps})$ to first order in $\sqrt{\eps}$, we find that in the limit $\eps \downarrow 0$, the rescaled fluctuation $\eta_t = \lim_{\eps \downarrow 0} \eta_t^{\eps}$ satisfies
\begin{equation}
    d \eta_t = \nabla b(t, X_t)\, \eta_t\, dt + \sqrt{2}\,\sigma(t, X_t)\,dW_t.
\end{equation}
Since this equation is linear in $\eta_t$ and driven by Gaussian noise, the variables $\eta_t$ are Gaussian. This tells us that $X_{t}^{\eps} = X_{t} + \sqrt{\eps}\,\eta_{t} + o(\sqrt{\eps})$. Large deviations go beyond this Gaussian description by characterizing fluctuations of order $\mathcal{O}(1)$ away from the deterministic path.

\subsubsection{Large Deviation Principle (LDP)}

To probe fluctuations beyond the CLT scale, consider the exponentially weighted expectation
\begin{equation}
    A^{\eps} = \mathbb{E}^{x_0}\left[e^{f(X_T^{\eps})/\eps} \right],
    \label{eq:A_definition}
\end{equation}
for some function $f:\mathbb{R}^d \to \mathbb{R}$. The factor $1/\eps$ in the exponent ensures that $A^{\eps}$ is dominated by realizations where $X_{t}^{\eps}$ deviates by $\mathcal{O}(1)$ from the deterministic trajectory---these are the \textit{large deviations} we seek to characterize. Our aim is to show that:
\begin{propbox}
\begin{equation}
    \lim_{\eps \downarrow 0} \eps \log A^{\eps} = \sup_{\{x_t\}_{t\in[0,T]}} \left[f(x_T) - \int_0^T L(t, x_t, \dot{x}_t)\, dt \right],
    \label{eq:A_eps_variational}
\end{equation}
\end{propbox}
where the supremum is over paths with $x_0$ fixed and the \textit{Lagrangian} is
\begin{equation}
    \label{eq:Lagrangian}
    L(t,x,\dot{x}) = \frac14(\dot{x} - b(t,x))\cdot D^{-1}(t,x) (\dot{x} - b(t,x)),
\end{equation} 
which measures the ``cost'' of deviating from the deterministic flow. We derive this  result systematically in the following subsections; a complementary heuristic derivation  based on path integrals is given in Appendix~\ref{appendix:path_integral}.

\begin{remark} 
A direct Monte Carlo estimate of $A^{\eps}$ is extremely inefficient for small $\eps$:  the CLT tells us that typical sampled paths remain within $\mathcal{O}(\sqrt{\eps})$ of the deterministic trajectory and contribute negligibly to the exponentially weighted average. The LDP~\eqref{eq:A_eps_variational} provides a principled solution. Once we identify the optimal path $x_t^*$ that achieves the supremum, we can design importance sampling schemes that preferentially generate paths near $x_t^*$, transforming an exponentially hard sampling problem into a deterministic optimization followed by efficient local sampling. We return to computational methods in Section~\ref{sec:numerical}.
\end{remark}

\begin{remark}[Connection to classical large deviation theory]
\label{rem:GE}
The result~\eqref{eq:A_eps_variational} is an instance of Varadhan's lemma~\cite{varadhan1966} applied to the Freidlin--Wentzell LDP for trajectories. The rate function for paths is the \textit{action} $S_T[x] = \int_0^T L(t, x_t, \dot{x}_t)\, dt$; Varadhan's lemma converts this to the variational formula for exponential expectations. The Hamiltonian $H(x,\theta)$ plays the role of a local cumulant generating function, with the Lagrangian obtained via Legendre transform---this is the infinitesimal version of the G\"artner--Ellis theorem~\cite{gartner1977,ellis1984}, which itself builds on the classical result by Cram\'er~\cite{cramer1938}. For Brownian motion specifically, the trajectory LDP is Schilder's theorem~\cite{schilder1966}; the extension to SDEs with drift is due to Freidlin and Wentzell~\cite{FreidlinWentzell2012}.
\end{remark}

\subsection{Derivation of the LDP: WKB analysis of the backward equation}
\label{sec:WKB_derivation}

To establish the LDP~\eqref{eq:A_eps_variational}, we must compute $A^{\eps} = \mathbb{E}^{x_0}[e^{f(X_T^{\eps})/\eps}]$ in the limit $\eps \downarrow 0$. The key tool is the \textit{backward Kolmogorov equation} (BKE), which governs conditional expectations of the stochastic process.

\paragraph{The generator and backward equation.} The \textit{infinitesimal generator} $L$ of a diffusion process~\eqref{eq:SDE_definition} gives a deterministic equation for the evolution of expectations. 
Let $\phi$ be an arbitrary compactly supported smooth test function, then its expectation evolves as,
\begin{equation}
    (L\phi)(t,x) = b(t,x) \cdot \nabla \phi(x) + \eps D(t,x) : \nabla \nabla \phi(x),
    \label{eq:generator_def}
\end{equation}
where $D:\nabla \nabla \phi = \sum_{i,j} D_{ij} \partial_{x_i} \partial_{x_j} \phi$. The derivation from the Euler--Maruyama discretization is given in Appendix~\ref{appendix:ITO}. 
We can, of course, use the generator to compute conditional expectations. The partial differential equation known as the backward Kolmogorov equation formalizes this: if $u^{\eps}(t,x)$ satisfies
\begin{equation}
    \partial_t u^{\eps} + Lu^{\eps} = 0, \qquad u^{\eps}(T,x) = F(x),
    \label{eq:BKE}
\end{equation}
then $u^{\eps}(0,x) = \mathbb{E}^{x}[F(X_T^{\eps})]$. The proof follows from It\^o's formula and the martingale property of stochastic integrals; see Appendix~\ref{appendix:ITO}. The BKE evolves \textit{backward} in time: given terminal data $F(x)$ at time $T$, the solution $u^{\eps}(t,x)$ at earlier times gives the expected value of $F(X_T^{\eps})$ conditioned on $X_t^{\eps} = x$.

\paragraph{WKB Ansatz}

The Wentzel-Kramers-Brillouin (WKB) ansatz takes its name from efforts to solve the Schr\"odinger equation in settings where $\hbar$ can be treated as a small parameter, though the approximation was first used by Liouville and Green~\cite{olver1974}.
The formal expansion yields accurate results when there is a separation of scales between the wavelength of solutions and the length scales of the environment.
The WKB approach relies on an asymptotic expansion in which the ansatz supposes that the leading order term in the solution can be expressed as $u^{\eps}(t,x) = e^{v^{\eps}(t,x)/\eps}$. We directly employ the ansatz in this form below.

\paragraph{Application to the LDP.} To compute $A^{\eps} = \mathbb{E}^{x_0}[e^{f(X_T^{\eps})/\eps}]$, we solve the BKE with terminal condition $u^{\eps}(T,x) = e^{f(x)/\eps}$, so that $A^{\eps} = u^{\eps}(0, x_0)$. The exponential form of the terminal condition suggests the WKB ansatz
\begin{equation}
    u^{\eps}(t,x) = e^{v^{\eps}(t,x)/\eps}.
\end{equation}
Substituting into the BKE $\partial_t u^{\eps} + b \cdot \nabla u^{\eps} + \eps D:\nabla\nabla u^{\eps} = 0$ and dividing by $u^{\eps}/\eps$:
\begin{equation}
    \partial_t v^{\eps} + b \cdot \nabla v^{\eps} + \eps D: \nabla \nabla v^{\eps} + \nabla v^{\eps} \cdot D \nabla v^{\eps}= 0.
\end{equation}
In the limit $\eps \downarrow 0$, the viscous term $\frac{\eps}{2}D: \nabla \nabla v^{\eps}$ becomes negligible, yielding the \textit{Hamilton--Jacobi equation}
\begin{propbox}
\begin{equation}
    \partial_t v + H(t,x, \nabla v) = 0, \qquad v(T,x) = f(x),
    \label{eq:HJE}
\end{equation}
\end{propbox}
with Hamiltonian
\begin{propbox}
\begin{equation}
    H(t,x,\theta) = b(t,x) \cdot \theta + \theta\cdot D(t,x)\, \theta.
    \label{eq:Hamiltonian}
\end{equation}
\end{propbox}
Since $v=\lim_{\eps\downarrow 0} v^\eps$ and $v^\eps(0,x_0) = \eps \log A^\eps$ by construction, we obtain the large deviation principle:
\begin{propbox}
\begin{equation}
    \lim_{\eps\downarrow 0 } \eps \log A^\eps = v(0,x_0),
    \label{eq:LDP:HJE}
\end{equation}
\end{propbox}
where $v(0,x_0)$ is determined by solving the Hamilton--Jacobi equation~\eqref{eq:HJE}.

\begin{remark}
The formal WKB expansion is valid when characteristics of the Hamilton--Jacobi equation~\eqref{eq:HJE} do not cross. When they do cross---corresponding to caustics in the variational problem---the classical solution ceases to exist and $\nabla v$ becomes multi-valued. The physically relevant continuation is then the \emph{viscosity solution}, which can be characterized as the limit $\lim_{\eps \to 0} v^{\eps}$ of solutions to the regularized equation. For our purposes, we work in regimes where characteristics remain well-separated; the general theory is developed in~\cite{FreidlinWentzell2012,flemingsoner2006,crandallishiilions1992}.
\end{remark}

\subsection{Method of characteristics and Hamilton's equations}

The Hamilton--Jacobi equation~\eqref{eq:HJE} is a first-order nonlinear PDE that can be solved by the method of characteristics. Along a characteristic curve $x_t$, we set $\theta_t = \nabla v(t, x_t)$, and the pair $(x_t, \theta_t)$ satisfies Hamilton's equations:
\begin{propbox}
\begin{equation}
    \label{eq:Hamilton_eqs}
    \begin{aligned}
    \dot{x}_t &= \partial_\theta H = b(t,x_t) + 2D(t,x_t)\,\theta_t, \\
    \dot{\theta}_t &= -\partial_x H = -[\nabla b(t,x_t)]^{\top} \theta_t - \nabla_x\bigl(\theta_t \cdot D(t,x_t)\,\theta_t\bigr),      
    \end{aligned}
\end{equation}
\end{propbox}
where $[\nabla_x(\theta \cdot D\,\theta)]_k = \sum_{ij}\theta_i (\partial_{x_k} D_{ij}) \theta_j$. Setting $\theta = 0$ recovers the deterministic dynamics: $\dot{x} = b(x)$. When $D$ is independent of $x$, the equation for $\dot{\theta}_t$ simplifies to $\dot{\theta}_t = -[\nabla b]^{\top} \theta_t$. The boundary conditions are
\begin{propbox}
\begin{equation}
    x_{t=0} = x_0, \qquad \theta_{t=T} = \nabla f(x_T).
    \label{eq:boundary_conditions}
\end{equation}
\end{propbox}
\vspace{0.5cm}
\begin{remark}[Boundary conditions for Hamilton's equations]
Unlike the standard initial value problem where both $x$ and $\theta$ are specified at $t=0$, the boundary conditions~\eqref{eq:boundary_conditions} constitute a two-point boundary value problem: position is fixed at the initial time, while momentum is determined at the final time by the terminal function $f$. This structure arises naturally from the LDP: the initial condition reflects the known starting point of the process, while the final momentum encodes the rare event we wish to characterize. Solving this boundary value problem generally requires iterative methods rather than simple forward integration; see Section~\ref{sec:numerical} for computational approaches.
\end{remark}

The path $x_t$ solving~\eqref{eq:Hamilton_eqs}--\eqref{eq:boundary_conditions} is called the \emph{instanton} (or optimal fluctuation path); it determines the most likely trajectory for the rare event specified by the boundary conditions. This path achieves the supremum in the variational principle~\eqref{eq:A_eps_variational}.

To extract the LDP, we compute the evolution of $v$ along characteristics. Using $\theta_t = \nabla v(t,x_t)$ and the Hamilton--Jacobi equation:
\begin{equation}
    \frac{d}{dt} v(t, x_t) = \partial_t v + \dot{x}_t \cdot \nabla v = -H + (b+ 2D\theta_t) \cdot \theta_t = \theta_t \cdot D\,\theta_t,
\end{equation}
where all the functions are evaluated at $(t,x_t)$.
Integrating from $0$ to $T$:
\begin{equation}
    v(T, x_T) - v(0, x_0) = \int_0^T \theta_t \cdot D(t,x_t)\,\theta_t\, dt.
\end{equation}
Since $v(T, x_T) = f(x_T)$ and $\lim_{\eps \downarrow 0} \eps \log A^{\eps} = v(0, x_0)$, we obtain the LDP:
\begin{propbox}
\begin{equation}
    \lim_{\eps \downarrow 0} \eps \log A^{\eps} = f(x_T) -  \int_0^T \theta_t \cdot D(t,x_t)\,\theta_t\, dt,
    \label{eq:LDP_result}
\end{equation}
\end{propbox}
where $(x_t, \theta_t)$ is the solution of Hamilton's equations~\eqref{eq:Hamilton_eqs} with boundary conditions~\eqref{eq:boundary_conditions}.

\subsection{Physical interpretation and Lagrangian formulation}

From $\dot{x}_t = b + 2 D\theta_t$, we see that $\theta_t = \frac12 D^{-1}(t,x_t) (\dot{x}_t - b(t,x_t))$ represents the ``force'' required to deviate from the deterministic drift. The optimal $\theta_t$ can be roughly interpreted as the most likely noise realization that produces the most likely path achieving a given rare event. The Hamiltonian formalism, combined with Laplace's principle, tells us that for small $\eps$, knowledge of this maximum-likelihood path captures the leading-order behavior of rare-event probabilities.

The LDP admits an equivalent Lagrangian formulation. The Lagrangian is obtained from the Hamiltonian via the Legendre transform:
\begin{equation}
    L(t,x, \dot{x}) = \sup_{\theta}\left[\dot{x} \cdot \theta - H(t,x, \theta)\right] = \frac{1}{4}(\dot{x} - b(t,x))\cdot D^{-1}(t,x)\, (\dot{x} - b(t,x)).
\end{equation}
This confirms~\eqref{eq:Lagrangian}. Using the involution property of the Legendre transform we also know that
\begin{equation}
    H(t,x, \theta) = \sup_{\dot{x}}\left[\dot{x} \cdot \theta - L(t,x, \dot{x})\right] =  \theta \cdot D(t,x)\,\theta
\end{equation}
Therefore the variational principle~\eqref{eq:A_eps_variational} can then be written as
\begin{equation}
\label{eq:LDP-Lagrangian}
    \lim_{\eps \downarrow 0} \eps \log A^{\eps} = \sup_{\{x_t\}_{t\in[0,T]}} \left\{ f(x_T) - \int_0^T L(t,x_t, \dot{x}_t)\, dt \right\},
\end{equation}
and the Euler--Lagrange equations reproduce Hamilton's equations~\eqref{eq:Hamilton_eqs} (see Appendix~\ref{appendix:euler_lagrange}).

More generally, the probability that $X_T^{\eps}$ lies in a set $B \subset \mathbb{R}^{d}$ satisfies
\begin{equation}
    \lim_{\eps \downarrow 0} \eps \log \mathbb{P}^{x_0}\left[X_T^{\eps} \in B\right] = -\inf_{x_T \in B} \int_0^T L(t,x_t, \dot{x}_t)\, dt,
\end{equation}
where the infimum is over optimal paths ending in $B$.


\begin{remark}[Connection to optimal control]
\label{rem:optimal-control}
The variational problem~\eqref{eq:LDP-Lagrangian} has an equivalent interpretation in the language of stochastic optimal control. Consider the controlled diffusion
\begin{equation}
    dX_t = b(t,X_t)\,dt + \sqrt{2\epsilon}\,\sigma(t,X_t)\,(dW_t + u_t\,dt),
\end{equation}
where $u_t$ is a control process. The cost of applying control $u$ is $\frac{1}{2}\int_0^T |u_t|^2\,dt$, and the Hamilton--Jacobi equation~\eqref{eq:HJE} is precisely the Hamilton--Jacobi--Bellman equation for minimizing this cost subject to a terminal reward $f(X_T)$. The optimal control is $u_t^* = \sigma^\top \nabla v$, which corresponds to the momentum $\theta_t = D\nabla v$ along the instanton. This equivalence underlies importance sampling schemes (Section~\ref{sec:importance-sampling}): biasing the dynamics by the optimal control produces trajectories that efficiently sample the rare event. For systematic development of the control-theoretic perspective, see~\cite{flemingsoner2006,DupuisEllis1997}; for modern neural-network approaches exploiting this connection via path-measure divergences and forward-backward SDEs, see~\cite{NueskenRichter2021}.
\end{remark}

\subsection{Example: Ornstein--Uhlenbeck process}

Consider the Ornstein--Uhlenbeck process with $\sigma = 1$:
\begin{equation}
    dX^{\eps}_{t} = -X^{\eps}_{t}\,dt + \sqrt{2\eps}\,dW_t, \quad X_{t=0} = 0,
\end{equation}
and the observable $A^{\eps} = \mathbb{E}^{0}[\exp(\lambda X_T^{\eps}/\eps)]$ for $\lambda \in \mathbb{R}$.

\paragraph{Exact solution.} This process can be solved exactly. Setting $Y_t = e^t X_t^{\eps}$ and applying It\^o's formula, we find $dY_t = \sqrt{\eps}\,e^t\,dW_t$, which integrates to
\begin{equation}
    X_{t}^{\eps} = \sqrt{2\eps}\,e^{-t}\int_{0}^t e^{s}\,dW_s.
\end{equation}
Thus $X_t^{\eps} \sim \mathcal{N}(0, \eps(1-e^{-2t}))$, and a direct Gaussian integration yields
\begin{equation}
    \lim_{\eps \downarrow 0} \eps \log A^{\eps} = \frac{\lambda^2}{2}(1-e^{-2T}).
    \label{eq:OU_exact}
\end{equation}

\paragraph{Hamiltonian approach.} With $b(x) = -x$, $D = 1$, and $f(x) = \lambda x$, Hamilton's equations~\eqref{eq:Hamilton_eqs} simplify to
\begin{equation}
    \dot{x}_t = -x_t + 2\theta_t, \qquad \dot{\theta}_t = \theta_t,
\end{equation}
with boundary conditions $x_0 = 0$ and $\theta_T = \lambda$. The second equation gives $\theta_t = \lambda e^{t-T}$. Substituting into the first and solving with $x_0 = 0$:
\begin{equation}
    x_{t} = \lambda e^{-T}(e^{t} - e^{-t}).
\end{equation}
Evaluating the LDP~\eqref{eq:LDP_result} with $D = 1$:
\begin{align}
    \lim_{\eps \downarrow 0} \eps \log A^{\eps} &= \lambda x_T -  \int_0^T |\theta_t|^2\, dt \nonumber \\
    &= \lambda^2(1 - e^{-2T}) - \frac{\lambda^2}{2}(1 - e^{-2T}) = \frac{\lambda^2}{2}(1-e^{-2T}),
\end{align}
in agreement with~\eqref{eq:OU_exact}. Note that the optimal path $x_t$ differs qualitatively from typical trajectories, which fluctuate around zero on scale $\sqrt{\eps}$.

\section{Large-volume limit: Markov jump processes}
\label{Sec.MJP}

We now extend the large deviation framework to Markov jump processes (MJPs), which describe systems evolving through discrete transitions between states. Let $X_t \in \Omega \subset \mathbb{R}^d$ (e.g., $\Omega = \mathbb{Z}^d$). The dynamics is specified by $R$ reaction channels with rates $\lambda_j(X) \geq 0$ and displacement vectors $\nu_j \in \mathbb{R}^d$ for $j = 1, \ldots, R$. In an infinitesimal time $\Delta t$:
\begin{equation}
    X_{t+\Delta t} \,|\, X_t = 
    \begin{cases}
        X_t + \nu_j, & \text{with probability } \lambda_j(X_t)\,\Delta t + o(\Delta t),\\
        X_t, & \text{with probability } 1 - \sum_{j=1}^R \lambda_j(X_t)\,\Delta t + o(\Delta t).
    \end{cases}
\end{equation}

The generator of this process, analogous to~\eqref{eq:generator_def}, is
\begin{equation}
    (L\phi)(x) = \lim_{\Delta t \downarrow 0}\frac{1}{\Delta t}\mathbb{E}^{X_t=x}[\phi(X_{t+\Delta t}) - \phi(X_t)] = \sum_{j=1}^{R}\lambda_j(x)[\phi(x+\nu_j) - \phi(x)].
    \label{eq:generator_MJP}
\end{equation}
See Appendix~\ref{appendix:MJP} for a derivation and further discussion of MJP fundamentals.

\paragraph{Interpretation and numerical simulation.} The definition above has a transparent interpretation in terms of the Gillespie algorithm~\cite{Gillespie1977}, which simulates MJPs exactly. Given the current state $X_t$, the process waits for an exponentially distributed time $\tau \sim \text{Exp}(\Lambda(X_t))$, where $\Lambda(X_t) = \sum_j \lambda_j(X_t)$ is the total transition rate. When the jump occurs, reaction $j$ is selected with probability $\lambda_j(X_t)/\Lambda(X_t)$, and the state updates to $X_{t+\tau} = X_t + \nu_j$. This cycle then repeats from the new state. The memoryless property of the exponential distribution ensures that the process is Markov: the future evolution depends only on the current state, not on how long the system has been in that state.

\subsection{Examples}

\paragraph{The Schl\"ogl model.} This paradigmatic model of chemical bistability~\cite{schlogl1972chemical} involves species $A$, $B$, and $S$ with reactions
\begin{align}
     S &\underset{k_1}{\stackrel{k_0}{\rightleftharpoons}} B,\\
     A + 2S &\underset{k_2}{\stackrel{k_3}{\rightleftharpoons}} 3S.
\end{align}
In a well-mixed volume $V$ with fixed concentrations $c_A$ and $c_B$, the state variable is $N \in \mathbb{N}$, the number of $S$ molecules. The transitions $N \to N \pm 1$ occur with rates
\begin{align}
    \lambda_{+}(N) &= k_1 V c_B + k_3 c_A N(N-1)/V, \\
    \lambda_{-}(N) &= k_0 N + k_2 N(N-1)(N-2)/V^2.
\end{align}
These rates are extensive (proportional to system size). Introducing the intensive variable $x = N/V$ and the small parameter $\eps = 1/V$, we obtain
\begin{align}
    \lambda_{+}^{\eps}(x) &= \frac{1}{\eps}\left[k_1 c_B + k_3 c_A\, x(x-\eps)\right], \label{eq:schlogl_rate_+} \\
    \lambda_{-}^{\eps}(x) &= \frac{1}{\eps}\left[k_0 x + k_2 x(x-\eps)(x-2\eps)\right]. \label{eq:schlogl_rate_-}
\end{align}
The limit $\eps \to 0$ corresponds to the thermodynamic limit at fixed density.

\paragraph{Discrete diffusion on a ring.} Consider $M$ compartments with $n_i$ particles in compartment $i$, where $\sum_i n_i = N_T$ is conserved. Particles hop between adjacent compartments at rate $\alpha n_i$:
\begin{align}
    n_i \to n_i - 1, \quad n_{i\pm 1} \to n_{i\pm 1} + 1, \quad \text{rate } \alpha n_i.
\end{align}
The displacement vectors are $\nu^{R}_i = -e_i + e_{i+1}$ and $\nu^{L}_i = -e_i + e_{i-1}$, where $e_i$ is the $i$th unit vector. With periodic boundary conditions and intensive variables $x_i = n_i/N_T$, the small parameter is $\eps = 1/N_T$.

\subsection{Rescaled generator and backward Kolmogorov equation}

The rates in both examples are extensive: proportional to system size $1/\eps$. This is generic for well-mixed reaction systems, where the number of reactive events per unit time scales with the volume. Passing to intensive variables, the rates take the form $\eps^{-1}[\lambda_j(x) + O(\eps)]$, where the leading-order rates $\lambda_j(x)$ depend only on concentrations. The rescaled generator is
\begin{equation}
    (L^{\eps} \phi)(x) = \sum_{j=1}^{R} \frac{1}{\eps}[\lambda_{j}(x)+O(\eps)][\phi(x + \eps \nu_{j}) - \phi(x)].
\end{equation}

For the Schl\"ogl model, expanding~\eqref{eq:schlogl_rate_+}--\eqref{eq:schlogl_rate_-}:
\begin{align}
    \lambda_{+}(x) &= k_1 c_B + k_3 c_A\, x^2, \\
    \lambda_{-}(x) &= k_0 x + k_2 x^3.
\end{align}
For discrete diffusion: $\lambda_{i,L}(x) = \lambda_{i,R}(x) = \alpha x_i$.

As in the diffusive case (Section~\ref{sec:WKB_derivation}), the generator determines the evolution of conditional expectations via the backward Kolmogorov equation: if $u^{\eps}(t,x)$ solves
\begin{equation}
    \partial_t u^{\eps} + L^{\eps} u^{\eps} = 0, \qquad u^{\eps}(T,x) = F(x),
    \label{eq:BKE_MJP}
\end{equation}
then $u^{\eps}(0,x) = \mathbb{E}^{x}[F(X_T^{\eps})]$. The proof is analogous to the SDE case; see Appendix~\ref{appendix:MJP}. We now use the BKE to derive the Law of Large Numbers and the Large Deviation Principle.

\subsection{Law of Large Numbers and the law of mass action}

Setting $u(t,x) = \lim_{\eps \downarrow 0} u^{\eps}(t,x)$ in the BKE~\eqref{eq:BKE_MJP} and expanding $\phi(x + \eps\nu_j) = \phi(x) + \eps\,\nu_j \cdot \nabla\phi + O(\eps^2)$, we obtain
\begin{equation}
    \partial_t u + \sum_{j=1}^{R} \lambda_{j}(x)\, \nu_j \cdot \nabla u = 0.
\end{equation}
The characteristics of this first-order PDE satisfy the \textit{law of mass action}:
\begin{equation}
    \dot{x}_t = \sum_{j=1}^{R} \lambda_{j}(x_t)\, \nu_j.
    \label{eq:mass_action}
\end{equation}
For the Schl\"ogl model: $\dot{x} = k_1 c_B + k_3 c_A\, x^2 - k_0 x - k_2 x^3$. For discrete diffusion: $\dot{x}_i = \alpha(x_{i+1} + x_{i-1} - 2x_i)$, which is the discrete diffusion equation.

The convergence of MJPs to deterministic dynamics in the large-volume limit is established rigorously in~\cite{kurtz1972,ethierkurtz1986}; see also~\cite{kipnislandim1999} for interacting particle systems.

\subsection{Derivation of the LDP: WKB analysis for jump processes}

To derive the LDP, we apply the WKB ansatz to the BKE~\eqref{eq:BKE_MJP}. Consider $u^{\eps}(t,x) = \mathbb{E}^{x}[e^{f(X_{T-t}^{\eps})/\eps}]$ with terminal condition $u^{\eps}(0,x) = e^{f(x)/\eps}$. Setting $u^{\eps} = e^{v^{\eps}/\eps}$ and substituting into the BKE, we must evaluate $L^{\eps} e^{v^{\eps}/\eps}$. Using the rescaled generator and expanding $v^{\eps}(x + \eps\nu_j)$ around $x$:
\begin{equation}
    v^{\eps}(x + \eps\nu_j) = v^{\eps}(x) + \eps\,\nu_j \cdot \nabla v^{\eps}(x) + O(\eps^2),
\end{equation}
we obtain
\begin{equation}
    e^{v^{\eps}(x + \eps\nu_j)/\eps} = e^{v^{\eps}(x)/\eps} \cdot e^{\nu_j \cdot \nabla v^{\eps}(x) + O(\eps)}.
\end{equation}
Therefore,
\begin{equation}
    L^{\eps} e^{v^{\eps}/\eps} = \sum_{j=1}^{R} \frac{1}{\eps}[\lambda_j(x) + O(\eps)] \left[e^{v^{\eps}(x + \eps\nu_j)/\eps} - e^{v^{\eps}(x)/\eps}\right] = \frac{1}{\eps} e^{v^{\eps}/\eps} \sum_{j=1}^{R} [\lambda_j(x) + O(\eps)] \left[e^{\nu_j \cdot \nabla v^{\eps}} - 1\right].
\end{equation}
Substituting into the BKE $\partial_t u^{\eps} + L^{\eps} u^{\eps} = 0$, dividing by $e^{v^{\eps}/\eps}/\eps$, and taking $\eps \to 0$, we obtain the Hamilton--Jacobi equation:
\begin{propbox}
\begin{equation}
    \partial_t v + H(x, \nabla v) = 0, \qquad v(T,x) = f(x).
    \label{eq:HJ_MJP}
\end{equation}
\end{propbox}
with Hamiltonian
\begin{propbox}
\begin{equation}
    H(x,\theta) = \sum_{j=1}^{R} \lambda_j(x) \left[ e^{\nu_j \cdot \theta} - 1 \right].
    \label{eq:H_MJP}
\end{equation}
\end{propbox}
The associated Hamilton's equations are
\begin{propbox}
\begin{equation}
    \begin{aligned}
    \dot{x}_t &= \partial_\theta H = \sum_{j=1}^{R} \lambda_j(x)\, e^{\nu_j \cdot \theta}\, \nu_j, \\
    \dot{\theta}_t &= -\partial_{x} H = -\sum_{j=1}^{R} \nabla \lambda_j(x) \left[ e^{\nu_j \cdot \theta} - 1 \right].
    \end{aligned}
    \label{eq:Hamilton_MJP}
\end{equation}
\end{propbox}
Setting $\theta = 0$ in the first equation recovers the law of mass action~\eqref{eq:mass_action}.

Note that the Hamiltonian~\eqref{eq:HJ_MJP} differs fundamentally from the diffusive case~\eqref{eq:Hamiltonian}. Expanding for small $\theta$: $H(x,\theta) \approx \sum_j \lambda_j(x)[\nu_j \cdot \theta + \frac{1}{2}(\nu_j \cdot \theta)^2]$, which resembles~\eqref{eq:Hamiltonian} but with state-dependent diffusion. The non-quadratic nature of~\eqref{eq:HJ_MJP} means that large deviations are not captured by a Gaussian approximation.

\begin{remark}
As in the diffusive case (Remark~\ref{rem:GE}), the Hamiltonian~\eqref{eq:H_MJP} plays the role of a cumulant generating function: for a single jump with displacement $\nu$ occurring at rate $\lambda$, the term $\lambda(e^{\nu \cdot \theta} - 1)$ is precisely the CGF of a Poisson process. The Lagrangian obtained via Legendre transform gives the rate function, following the G\"artner--Ellis structure.
\end{remark}

\begin{figure}[t]
    \centering
    \includegraphics[width=0.6\linewidth]{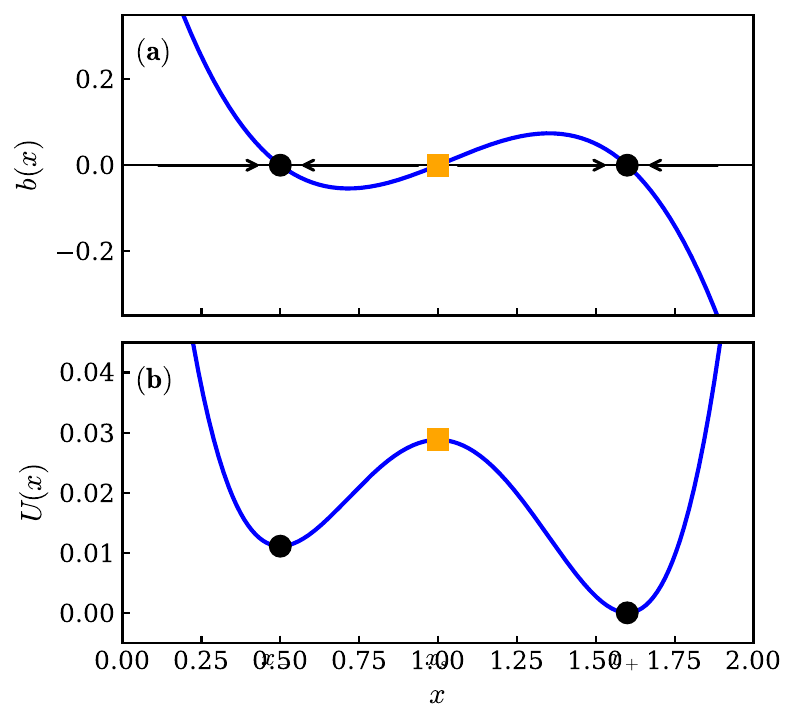}
    \caption{Deterministic dynamics of the Schl\"ogl model. Top: drift $b(x) = k_1 c_B + k_3 c_A\, x^2 - k_0 x - k_2 x^3$. Bottom: potential $U(x) = -\int b(x)\,dx$. Parameters: $k_1 c_B = 0.8$, $k_0 = 2.9$, $k_3 c_A = 3.1$, $k_2 = 1$. Black dots: stable fixed points; orange square: unstable fixed point. Arrows indicate the direction of relaxation.}
    \label{Fig.Schlogl}
\end{figure}

\paragraph{Schl\"ogl model.} With $\nu_{\pm} = \pm 1$:
\begin{equation}
    H(x,\theta) = (k_1 c_B + k_3 c_A\, x^2)(e^{\theta} - 1) + (k_0 x + k_2 x^3)(e^{-\theta} - 1).
    \label{eq:H_Schlogl}
\end{equation}
For appropriate parameter values, the deterministic dynamics $\dot{x} = \partial_\theta H|_{\theta=0}$ has two stable fixed points $x_{\pm}^*$ and one unstable fixed point $x_0^*$ (Fig.~\ref{Fig.Schlogl}).

\paragraph{Discrete diffusion.} The Hamiltonian is
\begin{equation}
    H(x,\theta) = \alpha \sum_{i=1}^{M} x_i \left[e^{\theta_{i+1}-\theta_i} + e^{\theta_{i-1}-\theta_i} - 2 \right].
    \label{eq:H_diffusion}
\end{equation}
Setting $\theta = 0$ recovers the discrete diffusion equation.

\paragraph{Continuum limit of discrete diffusion.} Dividing $[0,1]$ into $M$ boxes of width $\Delta s = 1/M$, we identify $s = i\,\Delta s$ and interpret $x_i \approx x(s)\,\Delta s$ as the fraction in compartment $i$, where $x(s)$ is a density satisfying $\int_0^1 x(s)\,ds = 1$. Similarly $\theta_i \to \theta(s)$. The physical diffusion coefficient is $D = \alpha(\Delta s)^2$, which we hold fixed as $\Delta s \to 0$. Expanding:
\begin{equation}
    e^{\theta_{i\pm 1} - \theta_i} - 1 = \pm \Delta s\, \partial_s\theta + \frac{(\Delta s)^2}{2}\left[\partial_s^2\theta + |\partial_s\theta|^2\right] + O((\Delta s)^3).
\end{equation}
Substituting into~\eqref{eq:H_diffusion} and using $\sum_i x_i[\ldots] = \sum_i x(s)\,\Delta s\,[\ldots] \to \int x(s)[\ldots]\,ds$:
\begin{equation}
    H[x,\theta] = D \int_{0}^1 ds\, x(s)\left[\partial_s^2 \theta + |\partial_s \theta|^2\right].
\end{equation}
The resulting Hamilton equations are PDEs for the fields $x(s,t)$ and $\theta(s,t)$.

\section{Long-time asymptotics I: The quasipotential}
\label{Sec.QP}

The quasipotential is a fundamental quantity for analyzing the long-time behavior of stochastic systems in the small-noise limit. It generalizes the notion of a potential landscape to nonequlibrium settings where detailed balance may not hold.

\subsection{Motivation: stationary distributions and effective dynamics in the small noise limit}

We consider autonomous systems---SDEs or MJPs with time-independent generators---in the small-noise or large-volume limit $\eps \to 0$. Our goal is to characterize:
\begin{itemize}
    \item The stationary distribution as $\eps \to 0$;
    \item When the deterministic dynamics from the LLN  has multiple stable fixed points, the effective dynamics over these metastable states.
\end{itemize}

The key object for both tasks is the \emph{quasipotential} $V(x, y)$, which measures the ``cost'' of a fluctuation-driven transition from $x$ to $y$. It is defined via a variational principle in Section~\ref{sec:variational_characterization} and generalizes the notion of a potential landscape to systems that may not satisfy detailed balance.

The quasipotential determines the stationary distribution as follows. For a system with a single stable fixed point $x_*$, the stationary distribution satisfies
\begin{equation}
    \rho_*(x) \asymp e^{-V(x)/\eps}.
    \label{eq:rho_QP}
\end{equation}
Here and throughout, $f \asymp g$ means $\lim_{\eps \to 0} \eps \log(f/g) = 0$, i.e., $f$ and $g$ have the same exponential scaling. In~\eqref{eq:rho_QP}, $V(x) = V(x_*, x)$. For systems with multiple stable fixed points $x_1, \ldots, x_n$, a single-argument function $V(x)$ can be constructed by patching together the functions $V(x_i, x)$ computed from each fixed point. The resulting $V(x)$ is the \emph{viscosity solution} of the stationary Hamilton--Jacobi equation (see Section~\ref{sec:QP_HJE}), and the stationary distribution still satisfies $\rho_*(x) \asymp e^{-V(x)/\eps}$.

The quasipotential also determines transition rates between metastable states. When the deterministic dynamics has multiple stable fixed points, the stochastic system exhibits \emph{metastability}: it spends long times near each fixed point, with rare noise-driven transitions between them. The transition rate from $x_i$ to an adjacent fixed point $x_j$ is
\begin{equation}
    k_{ij} \asymp e^{-V(x_i, x_j)/\eps}.
    \label{eq:Arrhenius_intro}
\end{equation}
For nonequlibrium systems, detailed balance is broken:
\begin{equation}
    \rho_*(x_i)\, k_{ij} \neq \rho_*(x_j)\, k_{ji}.
\end{equation}
This asymmetry---equivalently, $V(x_i, x_j) \neq V(x_j, x_i)$ in general---is a hallmark of systems that do not satisfy detailed balance and cannot be captured by the single-argument $V(x)$ alone.

\paragraph{Order of limits.} The quasipotential formalism arises from taking the limit $\eps \to 0$ first (yielding the Hamiltonian formalism of Sections~\ref{sec:LDP_trajectories} and~\ref{Sec.MJP}), followed by $T \to \infty$. These limits do not generally commute for multistable systems; see Section~\ref{sec:combined_limits} for a detailed discussion. The quasipotential characterizes the system's behavior on timescales $T \sim e^{c/\eps}$ for appropriate constants $c > 0$, during which transitions between metastable states become relevant. The systematic analysis of metastable transitions, including the hierarchy of timescales and the construction of cycles describing the sequence of most likely transitions, is developed in detail in~\cite{FreidlinWentzell2012}.

\subsection{The stationary Hamilton--Jacobi equation}
\label{sec:stationary_HJE}

To characterize the long-time behavior of the stochastic system, we seek the stationary distribution $\rho_*(x)$. For SDEs, this satisfies the stationary forward Kolmogorov equation
\begin{equation}
    0 = -\nabla \cdot (b\rho_*) + \eps \nabla \cdot (D \nabla \rho_*) = -\nabla \cdot J,
\end{equation}
where $J = b\rho_* - \eps D \nabla \rho_*$ is the probability current. Substituting the WKB ansatz $\rho_*(x) = w(x) e^{-V(x)/\eps}$ and collecting terms at leading order $O(\eps^{-1})$, we obtain the \emph{stationary Hamilton--Jacobi equation}:
\begin{propbox}
\begin{equation}
    H(x, \nabla V) = 0.
    \label{eq:stationary_HJE_QP}
\end{equation}
\end{propbox}
An analogous derivation starting from the stationary master equation yields the same structure for MJPs, with the Hamiltonian~\eqref{eq:HJ_MJP}.

A key feature of this equation is that it generically admits \emph{multiple solution branches}: the equation $H(x,\theta) = 0$ has several solutions $\theta(x)$ at each point $x$. In one dimension, SDEs with quadratic Hamiltonian~\eqref{eq:Hamiltonian} yield two branches; in higher dimensions and for MJPs, the structure can be more complex. 

One branch always corresponds to $\theta = 0$, which recovers the deterministic dynamics from the LLN. The other branches, with $\theta \neq 0$, correspond to fluctuation-driven motion against the deterministic flow.

These branches encode essential physics: each corresponds to optimal paths from a different stable fixed point, and the associated costs determine transition rates via~\eqref{eq:Arrhenius_intro}.

The HJE~\eqref{eq:stationary_HJE_QP} alone does not specify which branch to select at each point. To resolve this ambiguity, we define the quasipotential via a variational principle that explicitly specifies both endpoints.

\subsection{Variational definition of the quasipotential}
\label{sec:variational_characterization}

We now give a precise definition of the quasipotential that specifies both endpoints and thereby resolves the branch ambiguity.

Define the \emph{action} of a path $\{x_t\}_{t\in[0,T]}$ as
\begin{equation}
    S_T[x] = \int_0^T L(x_t, \dot{x}_t)\, dt,
\end{equation}
where $L(x, \dot{x}) = \sup_\theta[\dot{x} \cdot \theta - H(x,\theta)]$ is the Lagrangian obtained from the Hamiltonian via the Legendre transform. For SDEs with Hamiltonian~\eqref{eq:Hamiltonian}, this gives
\begin{equation}
    L(x, \dot{x}) = \frac{1}{4}(\dot{x} - b(x)) \cdot D^{-1}(x)\, (\dot{x} - b(x)),
\end{equation}
measuring the ``cost'' of deviating from the deterministic flow. For MJPs, the Lagrangian is obtained from~\eqref{eq:HJ_MJP} and has a more complex, non-quadratic form.

The minimum cost to travel from $x_a$ to $x_b$ in time $T$ is
\begin{equation}
    V_T(x_a, x_b) = \inf_{\substack{(x_t)_{t\in[0,T]}:\\ x_0=x_a, x_T=x_b}} S_T[x] \geq 0.
    \label{eq:V_T}
\end{equation}
The non-negativity follows from convexity of the Hamiltonian: $\dot{x} \cdot \theta - H(x,\theta) \geq -H(x,0) = 0$, with equality if and only if $\theta = 0$ (i.e., along the deterministic flow).

The \emph{quasipotential} is defined as the minimum cost over all possible transit times:
\begin{propbox}
\begin{equation}
    V(x, y) = \inf_{T > 0} V_T(x, y).
    \label{eq:quasipotential_var}
\end{equation}
\end{propbox}

This definition is unambiguous: $V(x, y)$ is the cost of the optimal fluctuation-driven path from $x$ to $y$, minimized over all possible durations.

\paragraph{Basic properties.} The quasipotential satisfies:
\begin{enumerate}
    \item \emph{Non-negativity}: $V(x, y) \geq 0$, with equality if and only if $y$ can be reached from $x$ along the deterministic flow.
    \item \emph{Identity}: $V(x, x) = 0$.
    \item \emph{Triangle inequality}: $V(x, z) \leq V(x, y) + V(y, z)$.
    \item \emph{Asymmetry} (generically): $V(x, y) \neq V(y, x)$ for nonequlibrium systems.
\end{enumerate}

\subsection{Optimal paths and the zero-energy surface}
\label{sec:QP_HJE}

We now show that for transitions between stable fixed points, the variational definition~\eqref{eq:quasipotential_var} recovers the stationary Hamilton--Jacobi equation~\eqref{eq:stationary_HJE_QP} along optimal paths.

\paragraph{Optimal paths between fixed points lie on $H = 0$.} Along any trajectory satisfying Hamilton's equations, the energy $E = H(x_t, \theta_t)$ is conserved. For a trajectory with constant energy $E \geq 0$, the action becomes
\begin{equation}
    S_T = \int_0^T \left[\dot{x}_t \cdot \theta_t - H\right] dt = \int_\gamma \theta \cdot dx - ET,
\end{equation}
where $\int_\gamma \theta \cdot dx$ is a line integral depending only on the geometric path $\gamma$, not on how it is traversed in time.

For paths connecting stable fixed points (where $\theta \to 0$ at the endpoints), the transit time $T \to \infty$ for $E = 0$, since the velocity $\dot{x} = \partial_\theta H \to 0$ as $\theta \to 0$. For $E > 0$, the term $-ET$ penalizes long transit times, but paths with $E > 0$ yield larger total action than the $E = 0$ paths; see~\cite{FreidlinWentzell2012} for a rigorous treatment.

Thus, for transitions between stable fixed points, the infimum in~\eqref{eq:quasipotential_var} is achieved along paths on the \emph{zero-energy surface} $H(x,\theta) = 0$:
\begin{equation}
    V(x_i, x_j) = \int_\gamma \theta \cdot dx,
    \label{eq:QP_line_integral}
\end{equation}
where $\gamma$ is the optimal path from $x_i$ to $x_j$ with $H(x,\theta) = 0$ throughout.

\paragraph{Recovery of the Hamilton--Jacobi equation.} Along the optimal path from a stable fixed point $x_*$ to $x$, the momentum satisfies $\theta = \nabla_x V(x_*, x)$. The constraint $H(x,\theta) = 0$ along this path shows that the two-point quasipotential satisfies the stationary HJE~\eqref{eq:stationary_HJE_QP} with respect to its second argument. 

Note that the HJE specifies $\nabla V$, not $V$ itself. For a system with multiple stable fixed points, the single-argument function $V(x)$ appearing in the stationary distribution satisfies $\nabla V(x) = \nabla_x V(x_i, x)$ in the basin of $x_i$—this selects the appropriate branch. At basin boundaries, $\nabla V$ is discontinuous as the branch changes, but $V$ itself remains continuous. The global function $V(x)$ is the viscosity solution of~\eqref{eq:stationary_HJE_QP}.

\begin{remark}
Constructing a single-argument function $V(x)$ such that $\rho_*(x) \asymp e^{-V(x)/\eps}$ from the two-point quasipotentials $V(x_i, x)$ requires adding appropriate constants when patching solutions from different basins. These constants are determined by the condition that $\rho_*(x_i)$ must be the stationary distribution of the effective Markov jump process over fixed points with rates~\eqref{eq:Arrhenius_intro}. Recall that this quasipotential formalism arises from taking $\eps \to 0$ first, then $T \to \infty$---the opposite order from the tilted generator approach in Section~\ref{sec:LDPtime}. See~\cite{FreidlinWentzell2012,schuttler2024} for details.
\end{remark}

\subsection{Gradient systems and detailed balance}
\label{sec:gradient_systems}

Before analyzing general nonequlibrium systems, we examine the important special case of gradient systems, where the quasipotential takes a particularly simple form.

Consider the SDE
\begin{equation}
    dX_t^{\eps} = -D\nabla U(X_t^{\eps})\,dt + \sqrt{2\eps}\,\sigma\,dW_t,
\end{equation}
where $D = \sigma\sigma^T$ and $U(x)$ is a potential. This system satisfies \emph{detailed balance}: the probability current vanishes at stationarity, and the stationary distribution is the Boltzmann distribution
\begin{equation}
    \rho_*(x) \propto e^{-U(x)/\eps}.
\end{equation}

\paragraph{Quasipotential equals the potential.} For gradient systems, the two-point quasipotential is simply the potential difference:
\begin{equation}
    V(x, y) = U(y) - U(x).
    \label{eq:QP_gradient}
\end{equation}
To verify this, note that the Hamiltonian for this system is
\begin{equation}
    H(x,\theta) = -D\nabla U(x) \cdot \theta + \theta^T D\theta.
\end{equation}
The zero-energy surface $H(x,\theta) = 0$ has two branches: $\theta = 0$ (relaxation) and $\theta = \nabla U$ (fluctuation). Integrating along the fluctuation branch from $x$ to $y$:
\begin{equation}
    V(x, y) = \int_\gamma \theta \cdot dx = \int_\gamma \nabla U \cdot dx = U(y) - U(x).
\end{equation}

\paragraph{Time-reversal symmetry.} The simplicity of~\eqref{eq:QP_gradient} reflects a deeper structure: for gradient systems, the optimal fluctuation path from $x$ to $y$ is the \emph{time-reversal} of the deterministic relaxation path from $y$ to $x$. This follows from Hamilton's equations: on the fluctuation branch $\theta = \nabla U$, we have
\begin{equation}
    \dot{x} = \partial_\theta H = D\nabla U(x),
\end{equation}
which is the time-reversal of the relaxation dynamics $\dot{x} = -D\nabla U(x)$ (up to a rescaling of time). The fluctuation path retraces the relaxation path in reverse.

\paragraph{Symmetric barriers.} As a consequence, the quasipotential is symmetric. For two minima $x_\pm$ separated by a saddle $x_s$:
\begin{equation}
    V(x_-, x_s) = U(x_s) - U(x_-), \qquad V(x_+, x_s) = U(x_s) - U(x_+).
\end{equation}
The ratio of transition rates satisfies
\begin{equation}
    \frac{k_{-+}}{k_{+-}} = e^{-[V(x_-,x_+) - V(x_+,x_-)]/\eps} = e^{-[U(x_+) - U(x_-)]/\eps},
\end{equation}
consistent with detailed balance: $\rho_*(x_-)k_{-+} = \rho_*(x_+)k_{+-}$.

\paragraph{Example: double-well potential.} Consider the one-dimensional double-well potential 
\begin{equation}
    U(x) = \frac{1}{4}(x^2 - 1)^2,
\end{equation}
with minima at $x_\pm = \pm 1$ and a saddle at $x_s = 0$. The quasipotential barriers are
\begin{equation}
    V(x_-, x_s) = V(x_+, x_s) = U(0) - U(\pm 1) = \frac{1}{4},
\end{equation}
and the transition rates are equal: $k_{-+} = k_{+-} \asymp e^{-1/(4\eps)}$.

\subsection{Non-gradient systems}
\label{sec:non_gradient}

For systems that do not satisfy detailed balance, the quasipotential framework remains valid but its structure is richer. We now discuss the key features that distinguish non-gradient systems from the equilibrium case.

\paragraph{Broken time-reversal symmetry.} For gradient systems, the optimal fluctuation path from $x$ to $y$ is the time-reversal of the relaxation path from $y$ to $x$. This symmetry is broken for non-gradient systems: the forward and backward optimal paths are generically distinct, reflecting the nonequlibrium nature of the dynamics.

\paragraph{Asymmetric quasipotential.} As a consequence, the two-point quasipotential is asymmetric:
\begin{equation}
    V(x, y) \neq V(y, x).
\end{equation}
The cost of transitioning from $x$ to $y$ differs from the reverse. This asymmetry has physical consequences: the ratio of forward and backward transition rates,
\begin{equation}
    \frac{k_{ij}}{k_{ji}} = e^{-[V(x_i, x_j) - V(x_j, x_i)]/\eps},
\end{equation}
is not determined by an equilibrium free energy difference, but by the full nonequlibrium dynamics. Detailed balance is violated:
\begin{equation}
    \rho_*(x_i)\, k_{ij} \neq \rho_*(x_j)\, k_{ji}.
\end{equation}

\paragraph{Non-trivial branch structure.} For gradient systems, the global quasipotential $V(x) = U(x) - U(x_*)$ corresponds to a single branch of the HJE (the fluctuation branch $\theta = \nabla U$). For non-gradient systems, the fluctuation branches emanating from different stable fixed points are genuinely distinct and must be explicitly patched together at basin boundaries, as described in Section~\ref{sec:QP_HJE}.

\paragraph{Limiting Markov chain on attractors.} When the deterministic dynamics has multiple stable fixed points $\{x_1, \ldots, x_n\}$, the quasipotential induces a hierarchy of timescales. On times of order $e^{V_{ij}/\eps}$, where $V_{ij} = V(x_i, x_j)$ is the cost of transition from $x_i$ to $x_j$, transitions between basins become relevant. The long-time dynamics can be approximated by a Markov jump process on the discrete set of attractors, with transition rates given by the Arrhenius formula~\eqref{eq:Arrhenius_intro}. 

This reduced Markov chain inherits the asymmetry of the quasipotential: its transition matrix is not reversible unless the original system satisfies detailed balance. In particular, the stationary distribution of this chain determines the constants needed to patch the single-argument quasipotential $V(x)$, as discussed in Section~\ref{sec:QP_HJE}. The construction of this limiting Markov chain, and the identification of the hierarchy of metastable states through cycles of transitions, is developed systematically in~\cite{FreidlinWentzell2012,schuttler2024}.

\subsection{Example: the Schl\"ogl model}
\label{sec:schlogl_example}

The Schl\"ogl model introduced in Section~\ref{Sec.MJP} provides a concrete illustration of the quasipotential formalism for non-gradient systems.

\subsubsection{Scalar case}

Recall that the Hamiltonian for the Schl\"ogl model is
\begin{equation}
    H(x,\theta) = \lambda_+(x)(e^{\theta} - 1) + \lambda_-(x)(e^{-\theta} - 1),
\end{equation}
with $\lambda_+(x) = k_1 c_B + k_3 c_A\, x^2$ and $\lambda_-(x) = k_0 x + k_2 x^3$. The deterministic drift is 
\begin{equation}
    b(x) = \lambda_+(x) - \lambda_-(x) = k_1 c_B + k_3 c_A\, x^2 - k_0 x - k_2 x^3.
\end{equation} 
For appropriate parameter values, the deterministic dynamics has two stable fixed points $x_\pm^*$ and one unstable fixed point $x_0^*$ (Fig.~\ref{Fig.Schlogl}).

\paragraph{Graphical construction.} In one dimension, the quasipotential has a simple graphical interpretation in the $(x,\theta)$ phase plane. The zero-energy surface $H(x,\theta) = 0$ consists of curves connecting the fixed points. For the Schl\"ogl model, solving $H(x,\theta) = 0$ yields a transcendental equation with two branches: the relaxation branch $\theta = 0$ (lying along the $x$-axis) and a fluctuation branch $\theta(x) \neq 0$.

\begin{figure}
    \centering
    \includegraphics[width=.6\linewidth]{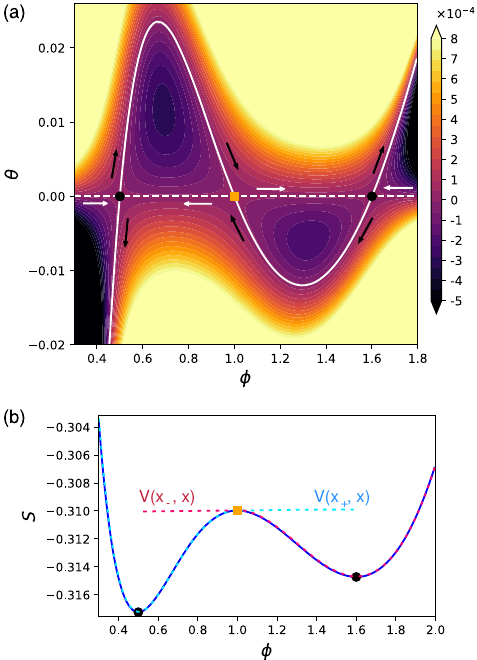}
    \caption{(a) Phase portrait for the bistable Schl\"ogl model showing the Hamiltonian contours $H(x, \theta)$. Black disks indicate stable fixed points $x_-$ (left) and $x_+$ (right); the orange square marks the unstable fixed point at $x_s$. White lines represent the zero-energy surface $H(x, \theta)=0$: the dashed horizontal line ($\theta=0$) corresponds to deterministic relaxation, while the solid curve ($\theta \neq 0$) represents optimal fluctuation paths. Arrows indicate the direction of motion along each branch. The quasipotential $V(x_-, x)$ for $x \leq x_s$ equals the area between the two branches from $x_-$ to $x$. (b) The global quasipotential $V(x)$, obtained by patching $V(x_-, x)$ and $V(x_+, x)$ with appropriate constants. Adapted from~\cite{zakine2023minimum}.}
    \label{Fig.PhSp-Pot}
\end{figure}

Fig.~\ref{Fig.PhSp-Pot}(a) shows this phase portrait. Starting from a stable fixed point $x_-$, the quasipotential $V(x_-, x)$ is computed by integrating along the fluctuation branch:
\begin{equation}
    V(x_-, x) = \int_{x_-}^{x} \theta(x')\, dx'.
    \label{eq:QP_integral_Schlogl}
\end{equation}
Geometrically, this integral equals the area between the two branches in the phase portrait. The global quasipotential, shown in Fig.~\ref{Fig.PhSp-Pot}(b), is obtained by computing $V(x_-, x)$ and $V(x_+, x)$ from each stable fixed point and patching them together with appropriate constants as described in Section~\ref{sec:QP_HJE}.

\paragraph{Quasipotential versus na\"ive potential.} One might attempt to construct a potential via
\begin{equation}
    U(x) = -\int^x b(x')\, dx' = -k_1 x - \frac{k_3}{3}x^3 + \frac{k_0}{2}x^2 + \frac{k_2}{4}x^4,
\end{equation}
and na\"ively expect the stationary distribution to be $\rho(x) \propto e^{-U(x)/\eps}$. However, this is incorrect: the Schl\"ogl model is a Markov jump process, not a gradient diffusion, and the Hamiltonian is not quadratic in $\theta$. The true quasipotential $V(x)$, obtained by solving $H(x, V'(x)) = 0$ and patching appropriately, differs from $U(x)$.

\begin{figure}[htp!]
    \centering
    \includegraphics[width=0.8\columnwidth]{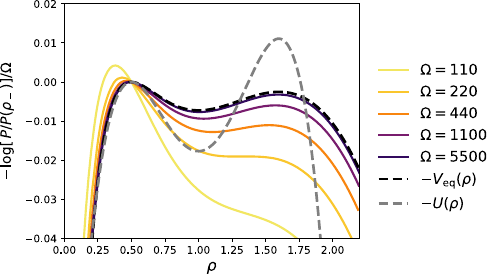}
    \caption{Rescaled logarithm of the stationary distribution for the Schl\"ogl model as the system size $\Omega = 1/\eps$ increases, showing convergence to the quasipotential $V(x)$ rather than the na\"ive potential $U(x)$. For the parameters shown, $x_-$ is the stable phase according to the quasipotential, while the na\"ive potential incorrectly predicts $x_+$ as stable. Parameters: $k_1 = 0.8$, $k_0 = 2.9$, $k_3 = 3.1$, $k_2 = 1$. Reproduced from~\cite{zakine2023minimum}.}
    \label{Fig:inv}
\end{figure}

Fig.~\ref{Fig:inv} shows that as the system size $\Omega = 1/\eps$ increases, the rescaled stationary distribution converges to $e^{-V(x)/\eps}$, not $e^{-U(x)/\eps}$. Strikingly, the relative stability of the two minima can be reversed: the na\"ive potential $U(x)$ may predict that $x_+$ is more stable, while the quasipotential correctly identifies $x_-$ as the globally stable state.

This discrepancy arises because the na\"ive potential $U(x)$ is constructed from the deterministic drift alone, ignoring the state-dependent fluctuation structure encoded in the full Hamiltonian. For gradient diffusions, the drift and fluctuations are related by detailed balance, and $U$ coincides with the quasipotential. For the Schl\"ogl model---and for non-gradient systems generally---this relationship breaks down.

\subsubsection{Spatially extended systems}

The Schl\"ogl model can be extended to a spatially distributed system~\cite{PhysRevE.86.040103} by dividing a domain into $L$ compartments. Within each compartment $i$, the Schl\"ogl reactions occur locally; particles also hop between adjacent compartments at rate $\alpha$. The Hamiltonian is
\begin{equation}
    H = H^R + H^D,
\end{equation}
where
\begin{align}
    H^R &= \sum_{i=1}^L \left[\lambda_{+}(x_i)(e^{\theta_i} - 1) + \lambda_{-}(x_i)(e^{-\theta_i} - 1)\right], \\
    H^D &= \alpha \sum_{i=1}^L x_i \left[e^{\theta_{i+1} - \theta_i} + e^{\theta_{i-1} - \theta_i} - 2\right].
\end{align}

\begin{figure}
    \centering
    \includegraphics[width=.6\linewidth]{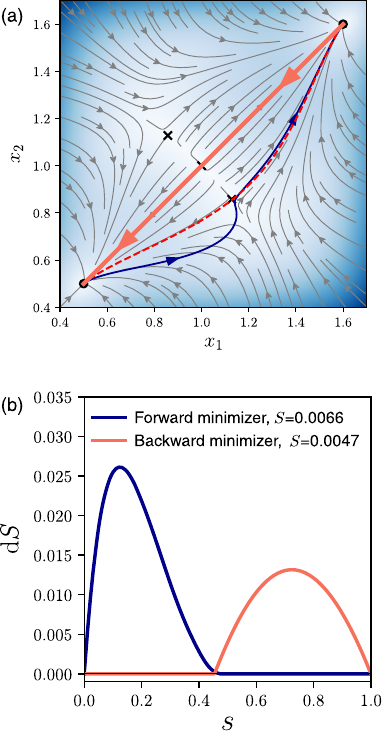}
    \caption{Optimal transition paths in the bistable Schl\"ogl model with two reactive compartments coupled by particle exchange at rate $\alpha$. (a) Flow lines of the deterministic dynamics with stable fixed points (black dots) and unstable fixed points (black crosses). The solid blue and orange lines show the minimum-action paths between metastable states; the dashed red line is the heteroclinic orbit. The divergence between these paths indicates nonequlibrium behavior. (b) Action accumulated along both minimum-action paths, demonstrating that the lower-left state is more stable under stochastic fluctuations. Reproduced from~\cite{zakine2023minimum}.}
    \label{fig:test2}
\end{figure}

For two compartments ($L = 2$), the phase portrait becomes two-dimensional (Fig.~\ref{fig:test2}). The forward and backward minimum-action paths connecting the metastable states are distinct, confirming nonequlibrium behavior. The asymmetry in the action along these paths determines which state is globally stable.

\begin{figure}[htp!]
    \centering
    \includegraphics[width=0.5\columnwidth]{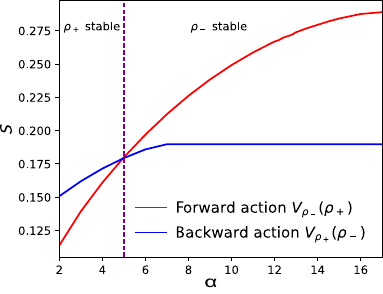}
    \caption{Phase diagram of the spatially extended Schl\"ogl model as a function of the hopping rate $\alpha$. For $\alpha < \alpha_c \approx 5.0$, the high-density phase $x_+$ is stable; for $\alpha > \alpha_c$, the low-density phase $x_-$ becomes stable. This nonequlibrium phase transition depends on the full stochastic dynamics, not just the deterministic flow. Parameters: $k_1 = 0.8$, $k_0 = 2.9$, $k_3 = 3.1$, $k_2 = 1$, $L = 40$. Reproduced from~\cite{zakine2023minimum}.}
    \label{Fig:Salpha}
\end{figure}

For larger systems, the diffusion rate $\alpha$ controls a nonequlibrium phase transition~\cite{PhysRevE.86.040103}. At a critical value $\alpha_c$, the relative stability of the high-density ($x_+$) and low-density ($x_-$) phases reverses (Fig.~\ref{Fig:Salpha}). This transition cannot be predicted from the deterministic dynamics alone---it emerges from the interplay of reactions and diffusion in the stochastic system.

\begin{figure}[htp!]
    \centering
    \includegraphics[width=0.8\columnwidth]{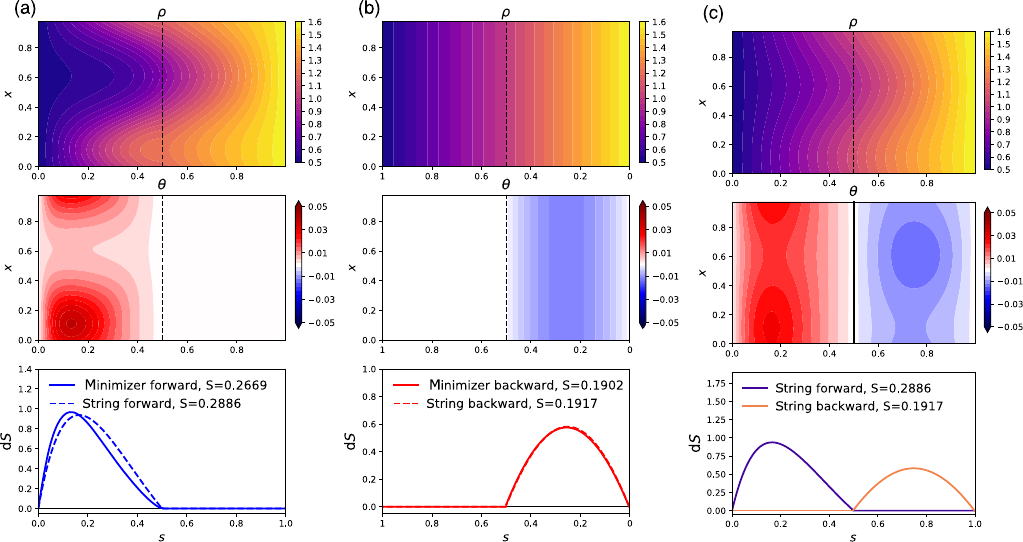}
    \caption{Optimal transition paths in the spatially extended Schl\"ogl model with $L = 40$ compartments. Top row: density profiles $\hat{\rho}(s, x)$ along the transition path (parameterized by $s$). Middle row: conjugate momentum $\hat{\theta}(s, x)$. Bottom row: Lagrangian along the trajectory. (a) Transition from $x_-$ to $x_+$. (b) Transition from $x_+$ to $x_-$. Note the pronounced asymmetry between forward and backward paths. (c) Heteroclinic orbit connecting $x_-$ and $x_+$. Parameters: $k_1 = 0.8$, $k_0 = 2.9$, $k_3 = 3.1$, $k_2 = 1$, $\alpha = 12$. Reproduced from~\cite{zakine2023minimum}.}
    \label{Fig:box40}
\end{figure}

Fig.~\ref{Fig:box40} shows the optimal transition paths for $L = 40$ compartments. The transitions $x_- \to x_+$ and $x_+ \to x_-$ proceed via qualitatively different mechanisms: one involves nucleation of a localized domain, while the other occurs more uniformly. This asymmetry is a hallmark of nonequlibrium phase coexistence.

\subsubsection{Pitfall: adding noise to the deterministic equations}

It is tempting to capture fluctuations by adding small noise to the deterministic equations from the LLN. For the continuum Schl\"ogl model, this procedure yields the Allen--Cahn/Ginzburg--Landau equation
\begin{equation}
    \partial_t x = \alpha \partial_s^2 x + \lambda_+(x) - \lambda_-(x),
\end{equation}
which is a gradient flow for the energy
\begin{equation}
    E[x] = \int_0^1 \left[\frac{\alpha}{2}|\partial_s x|^2 + U(x)\right] ds.
\end{equation}
If this were the correct description, the stationary distribution would be $\rho[x] \propto e^{-E[x]/\eps}$, implying that the relative stability of phases depends only on $U(x)$, not on $\alpha$.

However, the full stochastic analysis reveals an $\alpha$-dependent phase transition (Fig.~\ref{Fig:Salpha}). The na\"ive approach fails because it ignores the non-Gaussian fluctuation structure encoded in the full Hamiltonian~\eqref{eq:HJ_MJP}. This is a cautionary tale: \textit{adding noise to deterministic equations does not generally reproduce the correct large-deviation behavior}.

\section{Long-time asymptotics II: Finite noise and tilted generators}
\label{sec:LDPtime}

The preceding sections analyzed large deviations in the small-noise limit $\eps \to 0$, either at fixed time $T$ or in the limit as $T\to\infty$ taken \textit{a~posteriori}. We now consider the complementary regime: fixed noise strength with $T \to \infty$. This long-time limit characterizes time-averaged observables and connects to spectral properties of the generator.

\subsection{The tilted generator}

Consider the diffusion process
\begin{equation}
    dX_t = b(X_t)\,dt + \sqrt{2}\sigma(X_t)\,dW_t,
\end{equation}
and the exponentially weighted expectation
\begin{equation}
    I_T = \mathbb{E}\left[e^{\int_0^T f(X_t)\,dt}\right],
\end{equation}
where $f: \mathbb{R}^d \to \mathbb{R}$ is an observable. We seek the behavior of $I_T$ as $T \to \infty$.

Define the joint process $(X_t, A_t)$ where $\dot{A}_t = f(X_t)$ with $A_0 = 0$, so that $A_T = \int_0^T f(X_t)\,dt$ and $I_t = \mathbb{E}[e^{A_T}]$. Let $\rho(t,x, a)$ be the joint probability density function of $(X_t, A_t)$; it evolves according to the forward Kolmogorov equation
\begin{equation}
    \partial_t \rho = -\nabla \cdot (b\,\rho) + D:\nabla \nabla  \rho - f(x)\,\partial_a \rho,
\end{equation}
where $D = \sigma\sigma^\top$. The $e^a$-weighted marginal $\hat{\rho}(t,x) = \int_{-\infty}^{\infty} e^a\,\rho(t,x,a)\,da$ satisfies
\begin{equation}
    \partial_t \hat{\rho} = -\nabla \cdot (b\,\hat{\rho}) + D:\nabla \nabla\hat{\rho} + f(x)\,\hat{\rho},
    \label{eq:tilted_FPE}
\end{equation}
and $I_T = \int_{\mathbb{R}^d} \hat{\rho}(T,x)\,dx$.

Note that $\hat{\rho}$ does not preserve mass: $\int_{\mathbb{R}^d} \hat{\rho}(t,x)\,dx = \mathbb{E}[e^{A_t}] = I_t$, which grows or decays depending on the sign of $f$. The normalized density $\tilde{\rho}(t,x) = \hat{\rho}(t,x) / \int_{\mathbb{R}^d} \hat{\rho}(t,y)\,dy$ satisfies
\begin{equation}
    \partial_t \tilde{\rho} = -\nabla \cdot (b\,\tilde{\rho}) + D:\nabla \nabla \tilde{\rho} + f(x)\,\tilde{\rho} - \lambda(t)\,\tilde{\rho}_t,
    \label{eq:normalized_tilted}
\end{equation}
where $\lambda(t) = \int_{\mathbb{R}^d} f(x)\,\tilde{\rho}(t,x)\,dx$. A direct calculation shows that
\begin{equation}
 \int f(x)\,\tilde{\rho}(t,x)\,dx = \frac{\int_{\mathbb{R}^d}   f(x) \hat{\rho}(t,x)dx}
 {\int_{\mathbb{R}^d}  \hat{\rho}(t,x)dx} =  \frac{\mathbb{E}[f(X_t) e^{A_t} ] }
 {\mathbb{E}[ e^{A_t} ]}=  \frac{\mathbb{E}[f(X_t) e^{\int_0^t f(X_s)\,ds} ] }
 {\mathbb{E}[ e^{\int_0^t f(X_s)\,ds} ]}=\frac{d}{dt}\log I_t
\end{equation}
Thus $\lambda(t) = \frac{d}{dt}\log I_t$ and so $I_T = e^{\int_0^T \lambda_t\,dt}$.

In the long-time limit, assuming convergence to a stationary distribution $\tilde{\rho}(t,x) \to \tilde{\rho}_*(x)$, we have $\lambda(t) \to \lambda_*$ and
\begin{propbox}
\begin{equation}
    \lim_{T \to \infty} \frac{1}{T} \log I_T = \lambda_*.
    \label{eq:SCGF_limit}
\end{equation}
\end{propbox}
The quantity $\lambda_*$ is the principal eigenvalue (the eigenvalue with largest real part) of the \textit{tilted generator}
\begin{propbox}
\begin{equation}
    L_f \tilde{\rho} = -\nabla \cdot (b\,\tilde{\rho}) + D:\nabla\nabla \tilde{\rho} + f(x)\,\tilde{\rho} = \lambda_*\,\tilde{\rho}.
    \label{eq:tilted_eigenvalue}
\end{equation}
\end{propbox}

When $f(x) = s\,\phi(x)$, the quantity $\lambda_*(s)$ is the \textit{scaled cumulant generating function} (SCGF) of the time-averaged observable $\bar{\phi} = \frac{1}{T}\int_0^T \phi(X_t)\,dt$. The Legendre transform of $\lambda_* \equiv\lambda_*(s)$ gives the large-deviation rate function for $\bar{\phi}$. The spectral approach to large deviations via tilted generators originates in the work of Donsker and Varadhan~\cite{donskervaradhan1975}; see~\cite{chetritetouchette2015} for a modern treatment in the context of nonequilibrium processes.

\begin{remark}[G\"artner--Ellis theorem and its limitations]
The Legendre relationship between Hamiltonian and Lagrangian encountered in Sections~\ref{sec:LDP_trajectories} and~\ref{Sec.MJP} is the infinitesimal version of the G\"artner--Ellis theorem. Here, this structure appears globally: the SCGF $\lambda_*(s)$ is the cumulant generating function for the time-averaged observable, and the rate function is obtained via Legendre transform, $I(a) = \sup_s[sa - \lambda_*(s)]$~\cite{gartner1977,ellis1984}. This approach succeeds when $\lambda_*(s)$ is everywhere differentiable; the rate function is then convex. If $\lambda_*(s)$ has a non-differentiable point, the Legendre transform yields only a lower bound, and the true rate function may be non-convex. Such singularities arise when the system exhibits qualitatively different mechanisms for achieving different fluctuation magnitudes. See~\cite{touchette2009,dembo2009large,ellis1985} for comprehensive discussions.
\end{remark}

\subsection{Population dynamics interpretation}

Equation~\eqref{eq:normalized_tilted} admits an interpretation as a birth-death process on an ensemble of trajectories. The term $f(x)\tilde{\rho} - \lambda_t\tilde{\rho}$ acts as a source/sink: trajectories with $f(X_t) > \lambda(t)$ are ``born'' (replicated), while those with $f(X_t) < \lambda(t)$ ``die'' (are removed). The constraint $\lambda(t) = \int_{\mathbb {R}^d}  f(x)\,\tilde{\rho}(t,x)\,dx$ maintains normalization.

This interpretation underlies \textit{population dynamics} or \textit{cloning}  algorithms for computing large deviations~\cite{Cloning2,Cloning1}. The algorithm  maintains an ensemble of $N_c$ trajectories (``clones''), each evolving independently  according to the original dynamics. At regular intervals $\Delta t$, each clone at  position $x^{(i)}$ is assigned a weight $w_i = e^{f(x^{(i)})\Delta t}$. The ensemble  is then resampled with replacement according to these weights: clones in regions  where $f(x)$ is large are duplicated, while those where $f(x)$ is small are eliminated.  The total weight $W = \sum_i w_i$ fluctuates; its time-averaged growth rate converges  to $\lambda_*$ as $N_c \to \infty$ and $\Delta t \to 0$. The resampled ensemble approximates the tilted distribution $\tilde{\rho}_*$, enabling direct measurement of  observables conditioned on the rare event.

\subsection{Combined small-noise and long-time limits}
\label{sec:combined_limits}

Consider now the combined limit $T \to \infty$ and $\eps \to 0$ for the process
\begin{equation}
    dX_t = b(X_t)\,dt + \sqrt{2\eps}\,\sigma(X_t)\,dW_t,
\end{equation}
with
\begin{equation}
    I_T = \mathbb{E}\left[e^{\frac{1}{\eps}\int_0^T f(X_t)\,dt}\right].
\end{equation}
The tilted eigenvalue problem becomes
\begin{equation}
    -\nabla \cdot (b\,\tilde{\rho}_*) + \eps D:\nabla\nabla \tilde{\rho}_* + \frac{1}{\eps}f(x)\,\tilde{\rho}_* = \frac{\lambda_*}{\eps}\,\tilde{\rho}_*,
\end{equation}
and the SCGF is
\begin{equation}
    \lim_{\eps \to 0} \lim_{T \to \infty} \frac{\eps}{T} \log I_T = \lambda_*.
\end{equation}

The WKB ansatz $\tilde{\rho}_* = e^{-\phi_*/\eps}$ yields the stationary Hamilton--Jacobi equation
\begin{propbox}
\begin{equation}
    H(x, \nabla\phi_*) = \lambda_*,
    \label{eq:stationary_HJ}
\end{equation}
\end{propbox}
with Hamiltonian $H(x,\theta) = b(x) \cdot \theta +\theta\cdot D(x)\theta + f(x)$.

\begin{remark}[Non-commutativity of limits]
We stress that the limits $\eps \to 0$ and $T \to \infty$ of the scaled cumulant generating function
\begin{equation}
    \frac{\eps}{T} \log \mathbb{E}\left[e^{\frac{1}{\eps}\int_0^T f(X_t)\,dt}\right]
\end{equation}
do not generally commute. For systems with a unique attractor, the two limits typically agree. However, for multistable systems the order matters: taking $T \to \infty$ first allows the system to explore transitions between metastable states (with rates exponentially small in $1/\eps$), while taking $\eps \to 0$ first traps trajectories in a single basin, precluding such transitions. More precisely, the invariant measure at finite $\eps$ assigns weight to all attractors, but the $\eps \to 0$ limit of this measure depends on how $\eps$ and $1/T$ approach zero together~\cite{FreidlinWentzell2012}.
\end{remark}

\subsection{Separation of time scales}

A powerful application of the long-time limit is to systems with well-separated time scales~\cite{bouchet2016large}. Consider
\begin{align}
    dX_t^{\eps} &= a(X_t^{\eps}, Y_t^{\eps})\,dt, \\
    dY_t^{\eps} &= \frac{1}{\eps}b(X_t^{\eps}, Y_t^{\eps})\,dt + \frac{\sigma}{\sqrt{\eps}}\,dW_t,
\end{align}
where $X_t^\eps$ is a slow variable driven by the fast variable $Y_t^\eps$, and for simplicity we take $\sigma \in \mathbb{R}$ constant. The generator of this process is
\begin{equation}
    L^{\eps} = a(x,y) \cdot \nabla_x + \frac{1}{\eps}\left( b(x,y) \cdot \nabla_y + \frac{\sigma^2}{2} \Delta_y \right).
\end{equation}
We seek the LDP for $\mathbb{E}^{x,y}[e^{f(X_T^{\eps})/\eps}]$. The corresponding backward Kolmogorov equation is $\partial_t u^{\eps} + L^{\eps} u^{\eps} = 0$ with terminal condition $u^{\eps}(T,x,y) = e^{f(x)/\eps}$. Since the terminal condition depends only on $x$, we use the two-scale WKB ansatz
\begin{equation}
    u^{\eps}(t, x, y) = w(t, x, y)\, e^{u(t,x)/\eps},
\end{equation}
where the exponent $u$ depends only on the slow variable. Substituting into the BKE and collecting powers of $\eps$, we find at $O(\eps^{-1})$:
\begin{equation}
    a(x, y) \cdot (\nabla_x u)\, w + b(x, y) \cdot \nabla_y w + \frac{\sigma^2}{2}\Delta_y w = -(\partial_t u)\, w.
    \label{eq:fast_eigenvalue}
\end{equation}
In equation~\eqref{eq:fast_eigenvalue}, the quantities $x$ and $\nabla_x u$ enter as parameters, and the equation determines $w$ as a function of $y$. For a solution $w(y)$ to exist, the remaining quantity $\partial_t u$ cannot be arbitrary---there must be a solvability condition. The structure becomes transparent upon denoting $\theta =\nabla_x u$ and defining the tilted generator
\begin{equation}
    \mathcal{L}_\theta w \equiv a(x, y) \cdot \theta\, w + b(x, y) \cdot \nabla_y w + \frac{\sigma^2}{2}\Delta_y w,
    \label{eq:tilted_generator_fast}
\end{equation}
so that equation~\eqref{eq:fast_eigenvalue} becomes the eigenvalue problem $\mathcal{L}_\theta w = \lambda w$ with $\lambda = -\partial_t u$. That is, the solvability condition is that $\lambda$ must belong to the spectrum of $\mathcal{L}_\theta$. For ergodic fast dynamics, the relevant eigenvalue is the principal one $\lambda(x,\theta)$---the eigenvalue with largest real part---which is real, simple, and has a positive eigenfunction $w > 0$. This positivity ensures the WKB ansatz $u^{\eps} = w\, e^{u/\eps}$ remains well-defined. The solvability condition $-\partial_t u = \lambda(x, \nabla_x u)$ then yields the Hamilton--Jacobi equation for the slow variable:
\begin{propbox}
\begin{equation}
    \partial_t u + H_{\text{eff}}(x, \nabla_x u) = 0, \qquad H_{\text{eff}}(x, \theta) = \lambda(x, \theta).
    \label{eq:effective_H}
\end{equation}
\end{propbox}

The interpretation is that the fast process generates an effective Hamiltonian for the slow variable: the principal eigenvalue $\lambda(x,\theta)$ encodes a long-time average over the fast dynamics of $Y$, tilted by the momentum $\theta$. Large deviations of $X$ then follow from a small-noise analysis governed by this effective Hamiltonian.

\paragraph{Example.} Consider
\begin{align}
    dX_t^{\eps} &= |Y_t^{\eps}|^2 \,dt - V'(X_t^{\eps})\,dt, \\
    dY_t^{\eps} &= -\frac{\gamma(X_t^{\eps})}{\eps}\,Y_t^{\eps}\,dt + \frac{\sigma}{\sqrt{\eps}}\,dW_t,
\end{align}
where $V(x)$ is some potential, $\gamma(x) > 0$, and $\sigma \in \mathbb{R}$. For fixed $x$, the fast variable $Y$ is an Ornstein--Uhlenbeck process with stationary variance $\sigma^2/(2\gamma(x))$. 

For this system, $a(x,y) = y^2 - V'(x)$ and $b(x,y) = -\gamma(x)\, y$. The tilted generator~\eqref{eq:tilted_generator_fast} becomes
\begin{equation}
    \mathcal{L}_\theta w = (y^2 - V'(x))\, \theta\, w - \gamma(x)\, y\, \partial_y w + \frac{\sigma^2}{2} \partial_y^2 w.
\end{equation}
We seek solutions of $\mathcal{L}_\theta w = \lambda w$ using the Gaussian ansatz $w = e^{B y^2/2}$. Substituting and using $\partial_y w = By\, w$ and $\partial_y^2 w = (B + B^2 y^2) w$:
\begin{equation}
    (y^2 - V'(x))\, \theta - \gamma(x)\, B\, y^2 + \frac{\sigma^2}{2}(B + B^2 y^2) = \lambda.
\end{equation}
Collecting powers of $y^2$:
\begin{align}
    y^2: \quad & \theta - \gamma(x)\, B + \frac{\sigma^2}{2} B^2 = 0, \\
    y^0: \quad & -V'(x)\, \theta + \frac{\sigma^2}{2} B = \lambda.
\end{align}
The first equation determines $B$. Solving the quadratic:
\begin{equation}
    B = \frac{\gamma(x) \pm \sqrt{\gamma(x)^2 - 2\sigma^2 \theta}}{\sigma^2}.
\end{equation}
Since $H_{\text{eff}}$ is a cumulant generating function of the fast dynamics, it must be convex in $\theta$ by Jensen's inequality, which selects the minus sign. Substituting into the second equation yields the principal eigenvalue:
\begin{equation}
    H_{\text{eff}}(x, \theta) = \lambda(x,\theta) = -V'(x)\, \theta + \frac{1}{2}\left[\gamma(x) - \sqrt{\gamma(x)^2 - 2\sigma^2 \theta}\right].
\end{equation}
At the LLN level, the deterministic dynamics is $\dot{x} = \partial_\theta H_{\text{eff}}|_{\theta=0} = -V'(x) + \sigma^2/(2\gamma(x))$, corresponding to a ``thermal'' correction from the fast fluctuations. For the full LDP, the non-quadratic structure of $H_{\text{eff}}$ captures rare events where the slow variable is driven by atypical fluctuations of the fast variable.

\begin{figure}[t]
    \centering
    \includegraphics[scale=0.25]{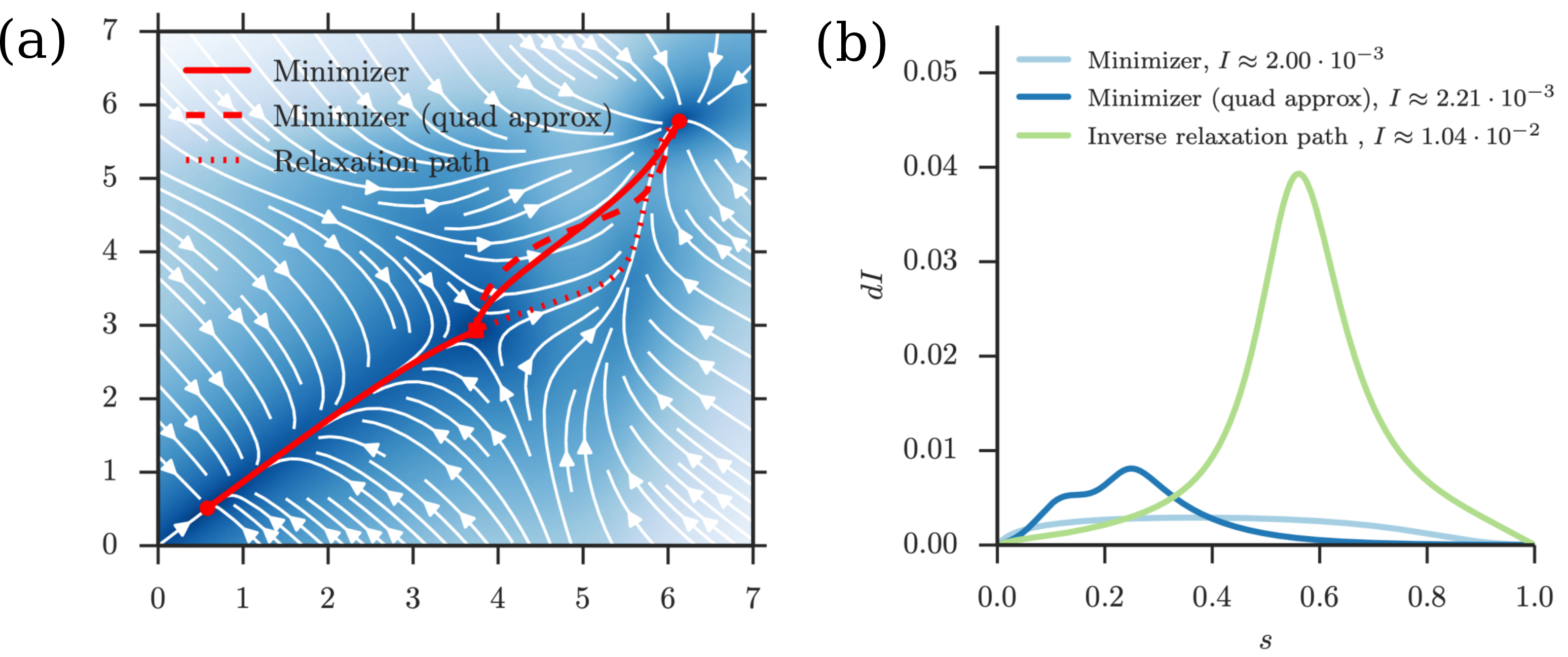}
    \caption{Optimal paths in a two-dimensional system with time-scale separation. (a) Phase portrait showing fixed points (dots) and the flow under the effective LLN dynamics. Solid curve: minimum-action path (instanton) between fixed points. Dotted curve: time-reversed relaxation path. Dashed curve: quadratic approximation to the Hamiltonian. (b) Cumulative action along each path. The instanton achieves approximately 20\% of the action required by the time-reversed path. Reproduced from~\cite{zakine2023minimum}.}
    \label{fig:time_separation_2d}
\end{figure}

Fig.~\ref{fig:time_separation_2d} illustrates a two-dimensional example where the optimal transition path differs significantly from a simple time-reversal of the relaxation path, with the noise-driven instanton achieving only about 20\% of the action of the naive reversed path.

\section{Computational methods}
\label{sec:numerical}

The variational problems arising from large deviation theory require numerical solution for all but the simplest models. This section surveys computational approaches, from classical minimum action methods to modern machine-learning techniques and sampling algorithms. Key references include~\cite{weinan2004mam} for the minimum action method,~\cite{vanden2008geometric} for its geometric version, and~\cite{PhysRevB.66.052301,ren2007simplified} for the string method.

\subsection{Minimum action methods}

The computation of instantons and quasipotentials reduces to finding paths that extremize the action functional~\eqref{eq:A_eps_variational}. Two equivalent formulations are available: the Lagrangian formulation, which directly minimizes the action
\begin{equation}
    S_T[x] = \int_0^T L(x_t, \dot{x}_t)\, dt,
    \label{eq:action_Lagrangian}
\end{equation}
and the Hamiltonian formulation, which seeks saddle points of
\begin{equation}
    S_T[x, \theta] = \int_0^T \left[\dot{x}_t \cdot \theta_t - H(x_t, \theta_t)\right] dt.
    \label{eq:action_Hamiltonian}
\end{equation}
Since $H$ is convex in $\theta$, this is a min-max problem: minimization over $x$ and maximization over $\theta$. Maximizing over $\theta$ at fixed $x$ recovers the Lagrangian via Legendre transform, $L(x, \dot{x}) = \max_\theta [\dot{x} \cdot \theta - H(x, \theta)]$. The Lagrangian formulation requires explicitly computing this transform, which assumes invertibility of the diffusion matrix $D$ for SDEs. The Hamiltonian formulation works directly with $H$, which is available from the generator without inversion, making it preferable when the diffusion is degenerate.

\subsubsection{The minimum action method (MAM)}

The minimum action method~\cite{weinan2004mam} discretizes the path $x_t$ as $\{x_0, x_1, \ldots, x_N\}$ at times $t_k = k\Delta t$, with $\Delta t = T/N$. The action is approximated by
\begin{equation}
    S_N = \sum_{k=0}^{N-1} L\left(\frac{x_k + x_{k+1}}{2}, \frac{x_{k+1} - x_k}{\Delta t}\right) \Delta t,
\end{equation}
with boundary conditions $x_0 = x_a$ and $x_N = x_b$ held fixed. The discretized action $S_N$ is minimized over the interior points $\{x_1, \ldots, x_{N-1}\}$ by gradient descent in a fictitious time $\tau$:
\begin{equation}
    \frac{\partial x_k}{\partial \tau} = -\frac{\delta S_N}{\delta x_k}.
\end{equation}
As $\tau \to \infty$, the path converges to a solution of the Euler--Lagrange equations. In practice, quasi-Newton methods such as L-BFGS significantly accelerate convergence over simple gradient descent.

When the diffusion matrix $D$ is not invertible and the Lagrangian cannot be computed explicitly, one works instead with the Hamiltonian formulation~\eqref{eq:action_Hamiltonian}; see e.g.~\cite{grafke2019numerical,zakine2023minimum}. Both $x$ and $\theta$ are discretized, and the saddle point is found by simultaneous gradient descent in $x$ and ascent in $\theta$. This gradient descent-ascent approach is developed systematically in Section~\ref{sec:GDA}.

The MAM is straightforward to implement for fixed transition time $T$, but cannot directly handle the limit $T \to \infty$ relevant for quasipotential calculations, since discretization points cluster near fixed points where the velocity vanishes. The geometric methods of Sections~\ref{sec:gMAM} and~\ref{sec:gGDA} address this limitation.

\subsubsection{Gradient descent-ascent for the min-max problem}
\label{sec:GDA}

A robust approach to solving the min-max problem~\eqref{eq:action_Hamiltonian} directly in the Hamiltonian formulation is gradient descent-ascent (GDA)~\cite{zakine2023minimum}. The problem is
\begin{equation}
    \inf_{x} S_T[x] = \min_x \max_\theta \int_0^T \left[\dot{x}_t \cdot \theta_t - H(x_t, \theta_t)\right] dt,
    \label{eq:minmax_IT}
\end{equation}
where the optimization is over trajectories $(x_t, \theta_t)_{t\in[0,T]}$ subject to boundary conditions $x_0 = x_a$ and $x_T = x_b$. The functional is strictly concave in $\theta$ (by convexity of $H$), so the min-max can be solved by alternating infinitesimal steps of steepest descent in $x$ and ascent in $\theta$. This leads to the GDA evolution equations in a fictitious optimization time $\tau$:
\begin{equation}
    \partial_\tau x = \alpha(\partial_t \theta + \partial_x H), \qquad
    \alpha \partial_\tau \theta = \partial_t x - \partial_\theta H,
    \label{eq:GDA_evolution}
\end{equation}
where $\alpha > 0$ is a parameter controlling the relative timescales. The fixed points of this system satisfy Hamilton's equations~\eqref{eq:Hamilton_eqs}, but the boundary conditions $x_{t=0} = x_a$ and $x_{t=T} = x_b$ must be imposed to select the instanton connecting the two states. This is again a two-point boundary value problem, but now with position specified at both endpoints rather than the mixed conditions~\eqref{eq:boundary_conditions} arising in the LDP; the momentum $\theta$ is determined implicitly. 

Since the system~\eqref{eq:GDA_evolution} is hyperbolic, and for such equations boundary conditions can only be imposed where characteristics enter the domain, we introduce the characteristic variables (not to be confused with the WKB ansatz of Section~\ref{sec:WKB_derivation})
\begin{equation}
    u = x + \alpha\theta, \qquad v = x - \alpha\theta,
\end{equation}
which satisfy
\begin{align}
    \partial_\tau u &= \partial_t u + f(u,v), \label{eq:u_evolution}\\
    \partial_\tau v &= -\partial_t v + g(u,v), \label{eq:v_evolution}
\end{align}
where
\begin{equation}
    f(u,v) = \alpha \partial_x H - \partial_\theta H, \qquad g(u,v) = \alpha \partial_x H + \partial_\theta H.
\end{equation}
The variable $u$ propagates backward in physical time $t$ (characteristics travel from $t=T$ to $t=0$), while $v$ propagates forward (characteristics travel from $t=0$ to $t=T$). Thus boundary conditions are naturally imposed at the upstream end of each characteristic: $u$ at $t=T$ and $v$ at $t=0$. Using $x = (u+v)/2$ and $\theta = (u-v)/(2\alpha)$, the constraints $x_{t=0} = x_a$ and $x_{t=T} = x_b$ become
\begin{equation}
    v(\tau, t=0) = 2x_a - u(\tau, t=0), \qquad u(\tau, t=T) = 2x_b - v(\tau, t=T).
\end{equation}
This suggests a natural algorithm based on Strang splitting: at each optimization step, first update~$u$ by propagating from $t=T$ to $t=0$, then update~$v$ by propagating from $t=0$ to $t=T$.

Algorithm~\ref{alg:MAM_GDA} gives the complete discretized scheme. The path is discretized on $M+1$ points with $\Delta t = T/M$, and an implicit upwind scheme ensures numerical stability.

\begin{algorithm}[H]
\caption{Action minimization by gradient descent-ascent}
\label{alg:MAM_GDA}
\begin{algorithmic}[1]
\REQUIRE Number of points $M$; initial path $\{x^0_i\}_{i \in \mathcal{I}}$ with $x^0_0 = x_a$ and $x^0_M = x_b$; functions $f(u,v)$ and $g(u,v)$; final time $T > 0$; step sizes $\Delta\tau > 0$, $\alpha > 0$.
\STATE \textbf{Initialize:} For all $i \in \mathcal{I} = \{0,\ldots,M\}$, set $\theta^0_i = 0$, $u^0_i = x^0_i + \alpha\theta^0_i$, $v^0_i = x^0_i - \alpha\theta^0_i$; set $\Delta t = T/M$.
\FOR{$n \geq 0$}
    \STATE \textbf{Update $u$ (backward sweep):} Solve $\{u^{n+1}_i\}_{i \in \mathcal{I}}$ sequentially from $i = M$ to $i = 0$:
    \begin{align*}
        u^{n+1}_M &= -v^n_M + 2x_b, \\
        \frac{u^{n+1}_i - u^n_i}{\Delta\tau} &= \frac{u^{n+1}_{i+1} - u^{n+1}_i}{\Delta t} + f(u^n_{i+1}, v^n_{i+1}), \quad i = M-1,\ldots,0.
    \end{align*}
    \STATE \textbf{Update $v$ (forward sweep):} Solve $\{v^{n+1}_i\}_{i \in \mathcal{I}}$ sequentially from $i = 0$ to $i = M$:
    \begin{align*}
        v^{n+1}_0 &= -u^{n+1}_0 + 2x_a, \\
        \frac{v^{n+1}_i - v^n_i}{\Delta\tau} &= -\frac{v^{n+1}_i - v^{n+1}_{i-1}}{\Delta t} + g(u^{n+1}_{i-1}, v^n_{i-1}), \quad i = 1,\ldots,M.
    \end{align*}
    \STATE \textbf{Recover fields:} Compute $x^{n+1}_i = \tfrac{1}{2}(u^{n+1}_i + v^{n+1}_i)$ and $\theta^{n+1}_i = \tfrac{1}{2\alpha}(u^{n+1}_i - v^{n+1}_i)$.
\ENDFOR
\end{algorithmic}
\end{algorithm}

The implicit upwind discretization is unconditionally stable and allows large time steps $\Delta\tau$, though the physical time step $\Delta t$ must still resolve the instanton dynamics. The reaction terms $f$ and $g$ are evaluated at upwind grid points, which improves stability for spatially extended fields.

\subsubsection{The geometric minimum action method (gMAM)}
\label{sec:gMAM}

For quasipotential calculations, optimal paths connect fixed points and have infinite duration ($T \to \infty$). The geometric minimum action method~\cite{vanden2008geometric} addresses this by exploiting the Maupertuis principle~\eqref{eq:QP_line_integral}: for paths on the zero-energy surface $H(x_t, \theta_t) = 0$, the action becomes the line integral
\begin{equation}
    S[x] = \int_\gamma \theta \cdot dx,
\end{equation}
which depends only on the geometric curve $\gamma$, not on its time parametrization. This reparametrization invariance allows treating infinite time intervals with finitely many discretization points.

The path is parametrized by normalized arc-length $s \in [0,1]$ with constant speed $\|x'(s)\|_D = L$, where the prime denotes differentiation with respect to $s$ and $\|v\|_D^2 = v \cdot D^{-1} v$ is the norm induced by the diffusion matrix. The constraint $H(x, \theta) = 0$ determines $\theta$ as a function of $(x, x')$:
\begin{equation}
    x' = \mu \nabla_\theta H(x, \theta),
    \label{eq:gMAM_constraint}
\end{equation}
where $\mu(s) = dt/ds$ is the time reparametrization factor. Enforcing the arc-length constraint gives
\begin{equation}
    \mu = \frac{\|x'\|_D^2}{\langle x', \nabla_\theta H \rangle_D},
\end{equation}
where $\langle u, v \rangle_D = u \cdot D^{-1} v$.

The gMAM algorithm iteratively updates the path:
\begin{enumerate}
    \item Solve $H(x^k, \theta^k) = 0$ for $\theta^k$ given $x^k$;
    \item Compute the descent direction $d^k = \theta^{k\prime} + \mu^k \nabla_x H(x^k, \theta^k)$;
    \item Update $x^{k+1} = x^k + \Delta\tau \, (\mu^k)^{-1} d^k$;
    \item Redistribute points to maintain equal arc-length spacing.
\end{enumerate}
The redistribution step prevents clustering of discretization points near fixed points where the velocity vanishes. The preconditioning by $\mu^{-1}$ ensures proper convergence where $\mu \to \infty$ as $\|b(x)\| \to 0$.

For further analytical and numerical approaches to computing the quasipotential in non-gradient systems, see~\cite{cameron2012finding,bouchet2016perturbative,schuttler2024}.

\subsubsection{Geometric GDA for quasipotential calculations}
\label{sec:gGDA}

The GDA approach of Section~\ref{sec:GDA} can be extended to compute the quasipotential by incorporating the geometric formulation~\cite{zakine2023minimum}. The problem becomes
\begin{equation}
    V_{x_a}(x_b) = \inf_{T>0} \inf_x S_T[x] = \min_{\lambda \geq 0} \min_{\hat x} \max_{\hat\theta} \int_0^1 \left[\hat x' \cdot \hat\theta - \lambda^{-1} H(\hat x, \hat\theta)\right] ds,
    \label{eq:geometric_minmax}
\end{equation}
where we reparametrize by $s \in [0,1]$ with $dt/ds = \lambda^{-1}(s)$, and $\hat x(s) = x(t(s))$, $\hat\theta(s) = \theta(t(s))$. By imposing $|\hat x'| = \text{const}$ (normalized arc-length), the discretization points remain uniformly distributed along the path regardless of physical time.

The GDA equations for the geometric formulation are
\begin{equation}
    \partial_\tau \hat x = \alpha\lambda \partial_s \hat\theta + \alpha \partial_{\hat x} H, \qquad
    \alpha \partial_\tau \hat\theta = \lambda \partial_s \hat x - \partial_{\hat\theta} H,
    \label{eq:geometric_GDA}
\end{equation}
where $\lambda$ is determined by the constraint that the path lies on the zero-energy surface $H = 0$ at convergence:
\begin{equation}
    \lambda = \frac{|\langle \partial_\theta H, \partial_s \hat x \rangle| + \sqrt{\max\bigl(0, \langle \partial_\theta H, \partial_s \hat x \rangle^2 - 4H|\partial_s \hat x|^2\bigr)}}{2|\partial_s \hat x|^2}.
    \label{eq:lambda_formula}
\end{equation}
At convergence, when $H(\hat x, \hat\theta) = 0$, the geometric Hamilton equations are satisfied:
\begin{equation}
    \lambda \hat x' = \partial_\theta H, \qquad \lambda \hat\theta' = -\partial_x H.
\end{equation}

Algorithm~\ref{alg:gMAM_GDA} gives the complete scheme for computing the quasipotential. The key modifications compared to Algorithm~\ref{alg:MAM_GDA} are the inclusion of the reparametrization coefficient $\lambda$ and the redistribution step.

\begin{algorithm}[H]
\caption{Geometric action minimization by gradient descent-ascent}
\label{alg:gMAM_GDA}
\begin{algorithmic}[1]
\REQUIRE Number of points $M$; stable fixed points $x_a$ and $x_b$ of the noiseless dynamics where $H(x_{a,b}, 0) = \partial_\theta H(x_{a,b}, 0) = 0$; initial path $\{\hat x^0_i\}_{i \in \mathcal{I}}$ with $\hat x^0_0 = x_a$, $\hat x^0_M = x_b$, and $|\hat x^0_{i+1} - \hat x^0_i|$ constant in $i$; functions $f(u,v)$, $g(u,v)$, and $\lambda(u,v)$; step sizes $\Delta\tau > 0$, $\alpha > 0$.
\STATE \textbf{Initialize:} For all $i \in \mathcal{I} = \{0,\ldots,M\}$, set $\hat\theta^0_i = 0$, $u^0_i = \hat x^0_i + \alpha\hat\theta^0_i$, $v^0_i = \hat x^0_i - \alpha\hat\theta^0_i$; set $\Delta s = 1/M$.
\FOR{$n \geq 0$}
    \STATE \textbf{Update $u$ (backward sweep):} Solve $\{u^{n+1}_i\}_{i \in \mathcal{I}}$ sequentially from $i = M$ to $i = 0$:
    \begin{align*}
        u^{n+1}_M &= -v^n_M + 2x_b, \\
        \frac{u^{n+1}_i - u^n_i}{\Delta\tau} &= \lambda_i(u^n, v^n) \frac{u^{n+1}_{i+1} - u^{n+1}_i}{\Delta s} + f(u^n_i, v^n_i), \quad i = M-1,\ldots,0.
    \end{align*}
    \STATE \textbf{Update $v$ (forward sweep):} Solve $\{v^{n+1}_i\}_{i \in \mathcal{I}}$ sequentially from $i = 0$ to $i = M$:
    \begin{align*}
        v^{n+1}_0 &= -u^{n+1}_0 + 2x_a, \\
        \frac{v^{n+1}_i - v^n_i}{\Delta\tau} &= -\lambda_i(u^{n+1}, v^n) \frac{v^{n+1}_i - v^{n+1}_{i-1}}{\Delta s} + g(u^{n+1}_i, v^n_i), \quad i = 1,\ldots,M.
    \end{align*}
    \STATE \textbf{Recover fields:} Compute $\bar x^{n+1}_i = \tfrac{1}{2}(u^{n+1}_i + v^{n+1}_i)$ and $\bar\theta^{n+1}_i = \tfrac{1}{2\alpha}(u^{n+1}_i - v^{n+1}_i)$.
    \STATE \textbf{Redistribute:} Interpolate $\{(\bar x^{n+1}_i, \bar\theta^{n+1}_i)\}_{i \in \mathcal{I}}$ onto a path $\{(\hat x^{n+1}_i, \hat\theta^{n+1}_i)\}_{i \in \mathcal{I}}$ such that $|\hat x^{n+1}_{i+1} - \hat x^{n+1}_i|$ is constant in $i$ (as in the string method).
    \STATE \textbf{Update characteristic variables:} Set $u^{n+1}_i = \hat x^{n+1}_i + \alpha\hat\theta^{n+1}_i$ and $v^{n+1}_i = \hat x^{n+1}_i - \alpha\hat\theta^{n+1}_i$.
\ENDFOR
\end{algorithmic}
\end{algorithm}

The redistribution step (line 6) prevents clustering of discretization points near fixed points where the velocity vanishes, while the coefficient $\lambda$ from~\eqref{eq:lambda_formula} automatically handles the transition from finite to infinite physical time by adjusting the speed along the path.

Both GDA algorithms require only evaluating the Hamiltonian $H(x, \theta)$ and its first derivatives $\partial_x H$ and $\partial_\theta H$---the same information needed to solve Hamilton's equations directly. This makes them applicable to systems with non-Gaussian noise (such as chemical reaction networks with Poisson fluctuations) where the Lagrangian is not available in closed form.

\subsubsection{The string method for gradient systems}

For gradient systems satisfying detailed balance (Section~\ref{sec:gradient_systems}), the optimal fluctuation path is the time-reversal of the deterministic relaxation path. As a result, instantons are heteroclinic orbits of the gradient flow $\dot{x} = \pm\nabla U(x)$, connecting two minima through the lowest saddle point. The string method~\cite{PhysRevB.66.052301,ren2007simplified} exploits this simplification.

In continuous time, the string evolves according to
\begin{equation}
    \partial_\tau x = -\nabla U(x) + \nu \hat{t},
    \label{eq:string_continuous}
\end{equation}
where $\tau$ is a fictitious relaxation time, $\hat{t} = x'/\|x'\|$ is the unit tangent vector, and $\nu$ is a Lagrange multiplier enforcing the parametrization constraint. The algorithm iteratively updates the path:
\begin{enumerate}
    \item Evolve each image: $\tilde{x}^{k+1} = x^k - \Delta\tau \, \nabla U(x^k)$;
    \item Redistribute points along the updated curve to maintain equal arc-length spacing.
\end{enumerate}
The redistribution step implicitly enforces the constraint without explicitly computing $\nu$. While restricted to gradient systems, the string method is computationally simpler than the gMAM and is widely used in chemical physics for computing minimum energy paths and transition states.

\subsection{Neural network approaches}

The methods described above require discretizing the path in time. For spatially extended systems governed by stochastic PDEs, the path must also be discretized in space, leading to high-dimensional optimization problems that are not only computationally costly but also technically challenging to implement: the underlying PDEs are typically stiff, requiring implicit time-stepping and careful treatment of boundary conditions. Neural network methods offer a compelling alternative by parametrizing the solution with neural networks and leveraging automatic differentiation for gradient computation~\cite{simonnet2023computing,zakine2023minimum}. 

\subsubsection{Deep MAM and deep gMAM}

Deep MAM and deep gMAM~\cite{zakine2023minimum,simonnet2023computing} adopt a physics-informed neural network (PINN) approach~\cite{raissi2019physics}: rather than discretizing the path on a mesh, they parameterize it by a neural network ansatz satisfying the boundary conditions. For a transition from state $x_a$ to state $x_b$, one constructs
\begin{equation}
    x_t(\omega) = (1-t/T)x_a + (t/T)x_b + (t/T)(1-t/T)\mathcal{N}_\omega(t),
\end{equation}
where $\mathcal{N}_\omega$ is a neural network with parameters $\omega$. This ansatz automatically satisfies $x_0(\omega) = x_a$ and $x_T(\omega) = x_b$. For field theories, $x_t$, $x_a$, $x_b$, and $\mathcal{N}_\omega$ are all functions of the spatial variable.

The action $S[x(\omega)]$ becomes a loss function over the network parameters, minimized by stochastic gradient descent. Algorithm~\ref{alg:deep_MAM} summarizes the procedure for finite transition time $T$.

\begin{algorithm}[H]
\caption{Deep MAM: Action minimization via neural network}
\label{alg:deep_MAM}
\begin{algorithmic}[1]
\REQUIRE Boundary states $x_a$, $x_b$; final time $T > 0$; neural network $\mathcal{N}_\omega$ with initial parameters $\omega^0$; number of collocation points $N_c$; learning rate $\eta$; number of iterations $N_{\text{iter}}$.
\STATE \textbf{Define path ansatz:} $x_t(\omega) = (1-t/T)x_a + (t/T)x_b + (t/T)(1-t/T)\mathcal{N}_\omega(t)$
\FOR{$n = 0, \ldots, N_{\text{iter}}-1$}
    \STATE \textbf{Sample collocation points:} Draw $\{t_j\}_{j=1}^{N_c}$ uniformly from $[0,T]$
    \STATE \textbf{Evaluate loss:} $\mathcal{L}(\omega^n) = \frac{T}{N_c} \sum_{j=1}^{N_c} L\bigl(x_{t_j}(\omega^n), \dot{x}_{t_j}(\omega^n)\bigr)$
    \STATE \textbf{Update parameters:} $\omega^{n+1} = \omega^n - \eta \nabla_\omega \mathcal{L}(\omega^n)$
\ENDFOR
\RETURN Optimized path $x_t(\omega^{N_{\text{iter}}})$
\end{algorithmic}
\end{algorithm}

The gradients $\nabla_\omega \mathcal{L}$ are computed via automatic differentiation, which also provides $\dot{x}_t = \partial_t x_t(\omega)$ at the collocation points. Space-time collocation points are randomly sampled at each optimization step, avoiding explicit grid discretization and the associated stability constraints.

For quasipotential calculations with $T \to \infty$, deep gMAM reparametrizes the path by arc-length $s \in [0,1]$ and minimizes the geometric action. This requires additional considerations: enforcing the arc-length constraint $|x'(s)| = \text{const}$ and the zero-energy condition $H(x,\theta) = 0$, typically via penalty terms or augmented Lagrangian methods. See~\cite{zakine2023minimum} for details.

\subsubsection{Physics-informed neural networks}

More broadly, physics-informed neural networks (PINNs)~\cite{raissi2019physics} can solve the Hamilton--Jacobi equation~\eqref{eq:HJE} or the stationary HJE~\eqref{eq:stationary_HJE_QP} directly. A neural network $v_\omega(t,x)$ is trained to minimize the residual
\begin{equation}
\mathcal{L}(\omega) = \int |\partial_t v_\omega + H(x, \nabla v_\omega)|^2 \, dx \, dt + \text{boundary terms},
\end{equation}
where the gradients are computed via automatic differentiation. This is useful for high-dimensional quasipotential calculations. For spatially extended systems where $x$ is itself a field, direct solution of the Hamilton--Jacobi equation becomes considerably more involved; score-based methods that learn the probability flow provide an alternative for computing probability currents and entropy production rates~\cite{boffi2023probability,nn_am}.

\subsection{Sampling methods for finite-noise, long-time limits}

The methods above compute instantons in the small-noise limit $\eps \to 0$. For the long-time limit at finite noise (Section~\ref{sec:LDPtime}), sampling-based approaches are needed. Genealogical algorithms~\cite{Cloning1,Cloning2} provide an efficient framework, and knowledge of the instanton can significantly improve their performance.

\subsubsection{Cloning algorithms}

The cloning algorithm maintains $N_c$ copies of the system, each evolving independently according to the original stochastic dynamics. At regular intervals $\Delta t$, each clone at position $x^{(i)}$ is assigned a weight $w_i = e^{f(x^{(i)})\Delta t}$, where $f(x)$ defines the tilted generator~\eqref{eq:tilted_eigenvalue}. The ensemble is resampled with replacement according to these weights: clones in regions where $f(x)$ is large are duplicated, while those where $f(x)$ is small are eliminated.

The total weight's growth rate converges to the principal eigenvalue $\lambda_*$ as $N_c \to \infty$:
\begin{equation}
    \lambda_* = \lim_{T \to \infty} \frac{1}{T} \log \mathbb{E}\left[\prod_{k} \bar{W}_k\right],
\end{equation}
where $\bar{W}_k$ is the average weight at the $k$th resampling step. The resampled ensemble samples the tilted distribution $\tilde{\rho}_*$, enabling measurement of observables conditioned on the rare event.

\subsubsection{Instanton-guided importance sampling}
\label{sec:importance-sampling}

Knowledge of the instanton can dramatically improve sampling efficiency through importance sampling. Given the instanton $(x^*_t, \theta^*_t)$, one biases the dynamics by adding a drift:
\begin{equation}
    dX_t = b(X_t)\,dt + 2D\theta^*_t\,dt + \sqrt{2\eps}\,\sigma\,dW_t.
\end{equation}
This keeps the original drift $b(X_t)$ responsive to the actual position while biasing toward the instanton. Samples are then reweighted using the Girsanov formula (Appendix~\ref{appendix:path_integral}) to correct for the modified dynamics, yielding unbiased estimates with reduced variance.

While this improves over direct sampling, using the instanton momentum $\theta^*_t$ as a time-dependent control does not achieve optimal variance reduction. Provably optimal schemes require adaptive controls that depend on the current position, as developed in~\cite{vanden2012rare}.

\subsection{Neural Optimal Control}
\label{sec:variational-control}

Similar approaches to instanton guidance can be employed within the Donsker-Varadhan framework, as well. 
For a dynamical observable of the form,
$A_T = \int_0^T f(X_t) \dd t + \int g(X_t) \circ \dd X_t,$
where $X_t$ evolves according to \eqref{eq:SDE_definition}, we can compute the scaled-cumulant generating function via an optimal importance sampling. 
As above, the scaled cumulant generating function 
\begin{equation}
    \lambda_*(s) = \lim_{T\to\infty} \frac1T \log \mathbb{E} \left[ e^{s A_T}\right].
\end{equation}
As shown in \cite{yan_learning_2022}, this problem can be reformulated as a variational optimal control problem in which one optimizes an optimal control drift, $f$ leading to a driven equation of motion,
\begin{equation}
    dX_t = b(X_t)\,dt + f(X_t)dt + \sqrt{2\eps}\,\sigma\,dW_t.
\end{equation}
The objective is obtained by constructing the optimal importance sampler and subsequently applying Jensen's inequality to obtain,
\begin{equation}
    \lambda_*(s) = \sup_{f} \lim_{T\to\infty} \left[ s\mathbb{E}_f [A_T] - T^{-1}\mathsf{D_{KL}(\mathbb{P}_f\|\mathbb{P})} \right],
    \label{eq:vardvldp}
\end{equation}
where $\mathsf{D_{KL}}$ denotes the Kullback--Leibler divergence between the path measures of the controlled and uncontrolled processes.
As developed in Refs.~\cite{yan_learning_2022, yan_physicsinformed_2022}, a neural network ansatz for $f$ yields a powerful and flexible ansatz for the control function. Directly maximizing \eqref{eq:vardvldp} via direct backpropagation or an adjoint method converges quickly and is robust even near dynamical phase transitions with appropriate sampling strategies.

\section*{Acknowledgments}

We thank the organizers of the Les Houches Summer School ``Large Deviations: Theory and Applications'' for the opportunity to present these lectures. W.D.P acknowlewdges funding from the Marie-Sk\l odowska-Curie grant No. 101104945.

\appendix

\section{Background on stochastic differential equations}
\label{appendix:ITO}

This appendix provides a self-contained introduction to the key results of stochastic calculus needed in these notes. We emphasize intuition and practical understanding over mathematical rigor; for complete treatments, we refer to~\cite{Oksendal2003,ELiVE2019,Gardiner2009}.

Throughout this appendix, we use uppercase letters $X_t$ and $\widetilde{X}_k$ to denote stochastic processes (continuous and discrete, respectively), and lowercase letters $x_t$ and $x_k$ for dummy integration variables representing paths.

Our starting point is the stochastic differential equation (SDE)
\begin{equation}
    dX_t = b(t, X_t)\,dt + \sigma(t,X_t)\,dW_t,
    \label{eq:app:SDE}
\end{equation}
where $X_t \in \mathbb{R}^d$, the drift $b:\mathbb{R}\times \mathbb{R}^d\to \mathbb{R}^d$ and volatility $\sigma:\mathbb{R}\times \mathbb{R}^d\to \mathbb{R}^{d\times d}$ are given functions, and $W_t$ is a $d$-dimensional Wiener process (Brownian motion). We denote the diffusion matrix by $D(t,x) = \sigma(t,x)\sigma(t,x)^\top$.

\subsection{Euler--Maruyama discretization}

The SDE~\eqref{eq:app:SDE} is defined as the continuum limit of a discrete-time process. To solve it numerically on a time interval $[0,T]$, divide time into $N$ steps of width $\Delta t = T/N$, with $t_k = k \Delta t$. The \emph{Euler--Maruyama scheme} defines the discrete approximation $\widetilde{X}_k \approx X_{t_k}$ by the recursion
\begin{propbox}
\begin{equation}
    \widetilde{X}_{k+1} = \widetilde{X}_{k} + b(t_k, \widetilde{X}_k) \Delta t + \sigma(t_k, \widetilde{X}_k)\,\sqrt{\Delta t}\,\xi_{k},
    \label{eq:EM}
\end{equation}
\end{propbox}
where $\xi_k \sim \mathcal{N}(0, \mathbb{1})$ are independent standard Gaussian random vectors. The noise term $\sqrt{\Delta t}\,\xi_k$ represents the Wiener increment
\begin{equation}
    \Delta W_k \equiv W_{t_{k+1}} - W_{t_k} \stackrel{d}{=} \sqrt{\Delta t}\,\xi_k,
\end{equation}
where $\stackrel{d}{=}$ denotes equality in distribution. This makes explicit that Brownian motion has increments of order $\sqrt{\Delta t}$, not $\Delta t$. The key properties are:
\begin{equation}
    \mathbb{E}[\Delta W_k] = 0, \qquad \mathbb{E}[\Delta W_k \Delta W_k^\top] = \Delta t\, \mathbb{1}, \qquad \mathbb{E}[\Delta W_k \Delta W_j^\top] = 0 \text{ for } k \neq j.
\end{equation}
As $\Delta t \to 0$, the discrete process converges to the continuous solution $X_t$, which is almost surely continuous but nowhere differentiable (H\"older continuous with exponent $1/2$). The $\sqrt{\Delta t}$ scaling of increments---rather than $\Delta t$---is the origin of the non-standard rules of stochastic calculus.

\subsection{The Markov property}

The Euler--Maruyama scheme~\eqref{eq:EM} makes manifest that the process is \emph{Markov}: given the current state $\widetilde{X}_k$, the next state $\widetilde{X}_{k+1}$ depends only on $\widetilde{X}_k$ and the fresh noise $\xi_k$, which is independent of all previous states and increments. The history of how the process reached $\widetilde{X}_k$ is irrelevant for predicting its future. 

This property carries over to the continuum limit: for any times $s < t$ and any measurable set $A\subset \mathbb{R}^d$,
\begin{equation}
    \mathbb{P}(X_t \in A \,|\, (X_r)_{ r \leq s}) = \mathbb{P}(X_t \in A \,|\, X_s).
\end{equation}
where $(X_r)_{ r \leq s}$ denotes the trajectory up to time $s$.
In words, conditioning on the entire trajectory up to time $s$ provides no more information about $X_t$ than conditioning on $X_s$ alone. The Markov property underlies all the results that follow, particularly the backward and forward Kolmogorov equations.

\subsection{The generator}

The \emph{generator} $\mathcal{L}$ of the diffusion process encodes the infinitesimal evolution of expected values. For any smooth test function $\phi: \mathbb{R}^d \to \mathbb{R}$, it is defined by
\begin{equation}
    (\mathcal{L}\phi)(x) = \lim_{\Delta t \downarrow 0} \frac{1}{\Delta t} \mathbb{E}^{X_t=x}\left[\phi(X_{t+\Delta t}) - \phi(x)\right].
    \label{eq:generator_def_app}
\end{equation}

\paragraph{Derivation from Euler--Maruyama.} Using~\eqref{eq:EM} and Taylor-expanding $\phi(\widetilde{X}_{k+1})$ around $\widetilde{X}_k = x$:
\begin{equation}
    \begin{aligned}
    \phi(\widetilde{X}_{k+1}) &= \phi\bigl(x + b\Delta t + \sigma\sqrt{\Delta t}\,\xi\bigr)  \\
    &= \phi(x) + \nabla\phi \cdot \bigl(b\Delta t + \sigma\sqrt{\Delta t}\,\xi\bigr)\\
    & + \frac{1}{2}\sum_{i,j} \partial_i\partial_j\phi\, (\sigma\sqrt{\Delta t}\,\xi)_i(\sigma\sqrt{\Delta t}\,\xi)_j + O(\Delta t^{3/2}).
\end{aligned}
\end{equation}
Since $\xi \sim \mathcal{N}(0,\mathbb{1})$, we have $\mathbb{E}[\xi] = 0$ and $\mathbb{E}[\xi_i\xi_j] = \delta_{ij}$. The cross term 
\[(\sigma\sqrt{\Delta t}\,\xi)_i(\sigma\sqrt{\Delta t}\,\xi)_j = \Delta t\,(\sigma\xi)_i(\sigma\xi)_j\] 
contributes at order $\Delta t$:
\begin{equation}
    \mathbb{E}[(\sigma\xi)_i(\sigma\xi)_j] = \sum_{k,l} \sigma_{ik}\sigma_{jl}\mathbb{E}[\xi_k\xi_l] = \sum_k \sigma_{ik}\sigma_{jk} = (\sigma\sigma^\top)_{ij} = D_{ij}.
\end{equation}
Thus:
\begin{equation}
    \mathbb{E}^{\widetilde{X}_{k}=x}[\phi(\widetilde{X}_{k+1}) - \phi(x)] = b \cdot \nabla\phi\, \Delta t + \frac{1}{2}D : \nabla\nabla\phi\, \Delta t + O(\Delta t^{3/2}).
\end{equation}
Dividing by $\Delta t$ and taking $\Delta t \to 0$ yields the generator:
\begin{propbox}
\begin{equation}
    (\mathcal{L}\phi)(x) = b(t,x) \cdot \nabla \phi(x) + \frac{1}{2} D(t,x) : \nabla \nabla \phi(x),
    \label{eq:generator}
\end{equation}
\end{propbox}
where $D : \nabla\nabla\phi = \sum_{i,j} D_{ij} \partial_i\partial_j\phi$. The second-order term arises because the noise enters at order $\sqrt{\Delta t}$, so its square contributes at order $\Delta t$.

\subsection{It\^o's formula}

It\^o's formula is the stochastic analogue of the chain rule. For a smooth function $\phi(t, x)$ evaluated along a trajectory $X_t$:
\begin{propbox}
\begin{equation}
    d\phi(t, X_t) = \left(\partial_t \phi + \mathcal{L}\phi\right)dt + \nabla\phi \cdot \sigma\,dW_t.
    \label{eq:Ito}
\end{equation}
\end{propbox}

\paragraph{Derivation:} From the Euler--Maruyama scheme, the increment over a small time step is $X_{t+\Delta t} - X_t \approx b\Delta t + \sigma\sqrt{\Delta t}\,\xi$ where $\xi \sim \mathcal{N}(0,\mathbb{1})$. Taylor-expanding $\phi(t+\Delta t, X_{t+\Delta t})$:
\begin{align}
    \phi(t+\Delta t, X_{t+\Delta t}) - \phi(t, X_t) &= \partial_t\phi\,\Delta t + \nabla\phi \cdot \bigl(b\Delta t + \sigma\sqrt{\Delta t}\,\xi\bigr) \nonumber \\
    &\quad + \frac{1}{2}\sum_{i,j}\partial_{x_i}\partial_{x_j}\phi\,(\sigma\sqrt{\Delta t}\,\xi)_i(\sigma\sqrt{\Delta t}\,\xi)_j + O(\Delta t^{3/2}).
\end{align}
The key observation is that
\begin{equation}
    (\sigma\sqrt{\Delta t}\,\xi)_i(\sigma\sqrt{\Delta t}\,\xi)_j = \Delta t\,(\sigma\xi)_i(\sigma\xi)_j = D_{ij}\Delta t + O(\Delta t),
\end{equation}
where the second equality holds in the mean-square sense: while $(\sigma\xi)_i(\sigma\xi)_j$ fluctuates, its deviation from $D_{ij}$ averages to zero and contributes only at higher order. Thus:
\begin{equation}
\begin{aligned}
    \phi(t+\Delta t, X_{t+\Delta t}) - \phi(t,X_t) &= \partial_t\phi\,\Delta t + b \cdot \nabla\phi\,\Delta t + \frac{1}{2}D:\nabla\nabla\phi\,\Delta t \\
    &+ \nabla\phi \cdot \sigma\sqrt{\Delta t}\,\xi + O(\Delta t^{3/2}).
\end{aligned}
\end{equation}
Taking $\Delta t \to 0$ and writing $\sqrt{\Delta t}\,\xi \to dW_t$ yields~\eqref{eq:Ito}. The ``extra'' term $\frac{1}{2}D:\nabla\nabla\phi$ compared to ordinary calculus arises because $(dW_t)^2 = dt$ rather than zero---this is the essence of It\^o calculus. In differential notation:
\begin{equation}
    (dW_t)^2 = dt, \qquad dW_t\,dt = 0, \qquad (dt)^2 = 0.
\end{equation}
These rules follow directly from the $\sqrt{\Delta t}$ scaling of Wiener increments.

\subsection{Martingale property and It\^o isometry}

Stochastic integrals $\int_0^T g(t,X_t) \cdot dW_t$ satisfy two fundamental properties. The first is the \emph{martingale property}:
\begin{propbox}
\begin{equation}
    \mathbb{E}\left[\textstyle\int_0^T g(t,X_t)  \cdot dW_t \right] = 0.
    \label{eq:martingale}
\end{equation}
\end{propbox}
The second is the \emph{It\^o isometry}:
\begin{propbox}
\begin{equation}
    \mathbb{E}\left[\left|\textstyle\int_0^T g (t,X_t) \cdot dW_t \right|^2\right] = \textstyle\int_0^T\mathbb{E}\left[ |g(t,X_t) |^2\right]\, dt.
    \label{eq:Ito_isometry}
\end{equation}
\end{propbox}
The term ``isometry'' reflects that stochastic integration preserves the $L^2$ norm: the variance of the stochastic integral equals the expected integral of $|g|^2$.

\paragraph{Derivation:} Discretizing $\int_0^T g(t,X_t) \cdot dW_t \approx \sum_k g_k \cdot \sqrt{\Delta t}\,\xi_k$ where $g_k = g(t_k, X_{t_k})$. Since $g_k$ depends only on the trajectory up to time $t_k$, it is independent of the current and future noise increments $\xi_j$ for $j \geq k$. In particular, $\mathbb{E}[\xi_k] = 0$ implies:
\begin{equation}
    \mathbb{E}\left[\sum_k g_k \cdot \sqrt{\Delta t}\,\xi_k\right] = \sqrt{\Delta t}\sum_k \mathbb{E}[g_k \cdot \xi_k] = \sqrt{\Delta t}\sum_k \mathbb{E}[g_k] \cdot \mathbb{E}[\xi_k] = 0.
\end{equation}

For the second moment:
\begin{equation}
    \mathbb{E}\left[\left|\sum_k g_k \cdot \sqrt{\Delta t}\,\xi_k\right|^2\right] = \Delta t \sum_{k,j} \mathbb{E}[(g_k \cdot \xi_k)(g_j \cdot \xi_j)].
\end{equation}
For $k < j$, the product $g_k \cdot \xi_k$ depends on the trajectory up to time $t_{k+1}$, and is therefore independent of $\xi_j$. Since $\mathbb{E}[\xi_j] = 0$:
\begin{equation}
    \mathbb{E}[(g_k \cdot \xi_k)(g_j \cdot \xi_j)] = \mathbb{E}[(g_k \cdot \xi_k)\, g_j] \cdot \mathbb{E}[\xi_j] = 0.
\end{equation}
By symmetry, the same holds for $k > j$. For $k = j$, using the independence of $g_k$ and $\xi_k$ together with $\mathbb{E}[\xi_{k,i}\, \xi_{k,l}] = \delta_{il}$:
\begin{equation}
    \mathbb{E}[(g_k \cdot \xi_k)^2] = \mathbb{E}\left[\sum_{i,l} g_{k,i}\, g_{k,l}\, \xi_{k,i}\, \xi_{k,l}\right] = \mathbb{E}\left[\sum_{i} g_{k,i}^2\right] = \mathbb{E}[|g_k|^2].
\end{equation}
Therefore:
\begin{equation}
    \mathbb{E}\left[\left|\sum_k g_k \sqrt{\Delta t}\,\xi_k\right|^2\right] = \Delta t \sum_k \mathbb{E}[|g_k|^2] \to\int_0^T \mathbb{E}\left[|g(t,X_t)|^2\right]\,dt \quad \text{as } \Delta t \to 0.
\end{equation}

\subsection{Backward Kolmogorov equation}

With It\^o's formula and the martingale property in hand, we can now derive the fundamental equations governing expectations over trajectories. The \emph{backward Kolmogorov equation} (BKE) governs the evolution of conditional expectations. Let $u(t,x) = \mathbb{E}^{X_t=x}[F(X_T)]$ for some terminal function $F$. Then $u$ satisfies
\begin{propbox}
\begin{equation}
    \partial_t u + \mathcal{L} u = 0, \qquad u(T,x) = F(x).
    \label{eq:BKE:b}
\end{equation}
\end{propbox}

\paragraph{Derivation:} By the Markov property and the tower rule for conditional expectations:
\begin{equation}
    u(t,x) = \mathbb{E}^{X_t=x}[F(X_T)] = \mathbb{E}^{X_t=x}\bigl[\mathbb{E}^{X_{t+\Delta t}}[F(X_T)]\bigr] = \mathbb{E}^{X_t=x}[u(t+\Delta t, X_{t+\Delta t})].
\end{equation}
Using the definition of the generator:
\begin{equation}
    u(t,x) = u(t+\Delta t, x) + \mathcal{L} u(t+\Delta t, x)\,\Delta t + O(\Delta t^2).
\end{equation}
Rearranging and taking $\Delta t \to 0$:
\begin{equation}
    \frac{u(t+\Delta t, x) - u(t,x)}{\Delta t} = -\mathcal{L} u + O(\Delta t) \quad \Rightarrow \quad \partial_t u = -\mathcal{L} u.
\end{equation}
The BKE evolves \emph{backward} in time: given terminal data $F(x)$ at time $T$, it determines $u(t,x)$ for $t < T$.

\begin{remark}
An alternative proof applies It\^o's formula directly to $u(t, X_t)$:
\begin{equation}
    du = [\partial_t u + \mathcal{L}u]\,dt + \nabla u \cdot \sigma\,dW_t.
\end{equation}
Integrating from $t$ to $T$ and taking expectations, the stochastic integral vanishes by the martingale property~\eqref{eq:martingale}:
\begin{equation}
    \mathbb{E}^{X_t=x}[u(T, X_T)] - u(t, x) = \mathbb{E}^{X_t=x}\left[\int_t^T [\partial_t u + \mathcal{L}u]\,ds\right].
\end{equation}
Since $u(T, X_T) = F(X_T)$ and $u(t,x) = \mathbb{E}^{X_t=x}[F(X_T)]$ by definition, the left-hand side vanishes. For this to hold for all $t$ and $x$, we must have $\partial_t u + \mathcal{L}u = 0$.
\end{remark}

\subsection{Forward Kolmogorov equation (Fokker--Planck)}

The \emph{probability density function} (PDF) $\rho(t,x)$ characterizes the distribution of $X_t$: for any measurable set $A \subseteq \mathbb{R}^d$,
\begin{equation}
    \mathbb{P}(X_t \in A) = \int_A \rho(t,x)\,dx.
\end{equation}
Equivalently, $\rho(t,x)\,dx$ is the probability of finding the process in an infinitesimal volume $dx$ around $x$ at time $t$. 

The \emph{forward Kolmogorov equation}, also called the \emph{Fokker--Planck equation}, governs the evolution of $\rho(t,x)$:
\begin{propbox}
\begin{equation}
    \partial_t \rho = \mathcal{L}^* \rho = -\nabla \cdot (b\,\rho) + \frac{1}{2}\nabla\nabla : (D\,\rho),
    \label{eq:FPE}
\end{equation}
\end{propbox}
where $\mathcal{L}^*$ is the adjoint of $\mathcal{L}$, and $\nabla\nabla:(D\rho) = \sum_{i,j}\partial_{x_i}\partial_{x_j}(D_{ij}\rho)$.

\paragraph{Derivation:} For any test function $\phi$, the expected value $\mathbb{E}[\phi(X_t)] = \int \phi(x)\rho(t,x)\,dx$ satisfies
\begin{equation}
    \frac{d}{dt}\int \phi\,\rho\,dx = \lim_{\Delta t \to 0} \frac{1}{\Delta t}\mathbb{E}[\phi(X_{t+\Delta t}) - \phi(X_t)] = \int (\mathcal{L}\phi)\,\rho\,dx.
\end{equation}
Integrating by parts (assuming boundary terms vanish):
\begin{equation}
    \int \phi\,\partial_t\rho\,dx = \int \phi\,(\mathcal{L}^*\rho)\,dx,
\end{equation}
which holds for all $\phi$, yielding~\eqref{eq:FPE}. The Fokker--Planck equation evolves \emph{forward} in time from an initial density $\rho(0,x)$.

\subsection{Path integral representation and Girsanov formula}
\label{appendix:path_integral}

The Euler--Maruyama discretization leads naturally to a path integral representation of expectations over trajectories. Consider for simplicity the SDE:
\begin{equation}
    dX_t = b(t,X_t)\,dt + dW_t,
\end{equation}
where for simplicity we take $\sigma = \mathbb{1}$. From~\eqref{eq:EM}, the update is
\begin{equation}
    \widetilde{X}_{k+1} =  \widetilde{X}_k + b(t_k, \widetilde{X}_k)\Delta t + \sqrt{\Delta t}\,\xi_k.
\end{equation}
The conditional distribution is
\begin{equation}
     \widetilde{X}_{k+1} \,|\,  \widetilde{X}_k \sim \mathcal{N}\left( \widetilde{X}_k + b(t_k, \widetilde{X}_k)\Delta t,\, \Delta t\, \mathbb{1} \right),
\end{equation}
with density proportional to $\exp\bigl(-\frac{1}{2\Delta t}|x_{k+1} - x_k - b(t_k,x_k)\Delta t|^2\bigr)$. The joint probability density function of a discretized path $(x_0, x_1, \ldots, x_K)$ is therefore
\begin{equation}
    P(x_1,\dots,x_K | x_0) \propto \exp\left(-\frac{1}{2\Delta t} \sum_{k=0}^{K-1}|x_{k+1} - x_{k} - b(t_k,x_k)\Delta t|^2\right).
\end{equation}

To interpret this, compare with pure Brownian motion ($b = 0$), for which 
\begin{equation}
    P_0(x_1,\dots,x_K | x_0) \propto \exp\left(-\frac{1}{2\Delta t} \sum_{k=0}^{K-1}|x_{k+1} - x_{k}|^2\right).
\end{equation}
The ratio of path probabilities is
\begin{align}
    \frac{P(x_1,\dots,x_K | x_0)}{P_0(x_1,\dots,x_K | x_0)} &= \exp\left(\sum_{k=0}^{K-1} b(t_k,x_k) \cdot (x_{k+1} - x_k) - \frac{1}{2}\sum_{k=0}^{K-1}|b(t_k,x_k)|^2\Delta t\right).
\end{align}

In the continuum limit $\Delta t \to 0$, recognizing $(x_{k+1} - x_k) \to dx_t$ formally and $\Delta t\sum_k  \to \int dt$,this becomes the \emph{Girsanov formula}:
\begin{propbox}
\begin{equation}
    \frac{\dd\mathbb{P}}{\dd\mathbb{P}_0}[(x_t)_{t \in [0,T]}] = \exp\left(\int_0^T b(t,x_t) \cdot \dd x_t - \frac{1}{2}\int_0^T |b(t,x_t)|^2\, \dd t\right),
    \label{eq:Girsanov}
\end{equation}
\end{propbox}
which gives the relative density of the path measure $\mathbb{P}$ (for the SDE with drift $b$) with respect to the Wiener measure $\mathbb{P}_0$. Here $(x_t)_{t \in [0,T]}$ denotes a continuous path starting at $x_0$, and $\dd x_t$ its infinitesimal increment. For a rigorous treatment, see~\cite{Oksendal2003,ELiVE2019}.

Using~\eqref{eq:Girsanov}, expectations under the SDE with drift $b$ can be expressed in terms of expectations under pure Brownian motion:
\begin{propbox}
\begin{equation}
\begin{aligned}
    &\mathbb{E}^{x_0}\left[F[(X_t)_{t\in[0,T]}]\right] = \mathbb{E}_0\left[F[(x_0 + W_t)_{t\in[0,T]}] M_T\right]\\
    \text{with} \quad &  M_T = \exp\left(\int_0^T b(t,x_0 + W_t) \cdot \dd W_t - \frac{1}{2}\int_0^T |b(t,x_0 + W_t)|^2\, \dd t\right),
\end{aligned}
    \label{eq:Girsanov_expectation}
\end{equation}
\end{propbox}
where $\mathbb{E}_0$ denotes expectation with respect to standard Brownian motion $W_t$ starting at $0$.
This exact formula has a formal  \emph{path integral representation}:
\begin{propbox}
\begin{equation}
    \mathbb{E}^{x_0}\left[F[(X_t)_{t\in[0,T]}]\right] \propto \int \mathcal{D}[x_t] \, F[(x_t)_{t\in[0,T]}] \exp\left(- \frac{1}{2} \int_{0}^T |\dot{x}_t - b(t,x_t)|^2\, dt\right),
    \label{eq:path_integral}
\end{equation}
\end{propbox}
where the integral is formally taken over continuous paths starting at $x_0$, and the Lagrangian $L(t,x,\dot{x}) = \frac{1}{2}|\dot{x} - b(t,x)|^2$ measures the ``cost'' of deviating from the deterministic flow. 

\subsection{Heuristic derivation of the LDP}

The path integral representation provides a direct route to the Large Deviation Principle. Consider now the SDE with small noise:
\begin{equation}
    dX_t^{\eps} = b(t,X_t^{\eps})\,dt + \sqrt{\eps}\,dW_t.
\end{equation}
Repeating the discretization argument, the path probability becomes
\begin{equation}
    P(x_1,\dots,x_K | x_0) \propto \exp\left(-\frac{1}{2\eps \Delta t} \sum_{k=0}^{K-1}|x_{k+1} - x_{k} - b(t_k,x_k)\Delta t|^2\right),
\end{equation}
leading to the path integral
\begin{equation}
    \mathbb{E}^{x_0}\left[e^{f(X_T^{\eps})/\eps}\right] \propto \int \mathcal{D}[x_t] \exp\left[ \frac{1}{\eps}\left(f(x_T) - \frac{1}{2} \int_{0}^T |\dot{x}_t - b(t,x_t)|^2\, dt\right)\right].
\end{equation}
As $\eps \to 0$, the factor $1/\eps$ in the exponent becomes large. By \emph{Laplace's method} (saddle-point approximation), the integral is dominated by paths that maximize the exponent:
\begin{propbox}
\begin{equation}
    \lim_{\eps \downarrow 0} \eps \log \mathbb{E}^{x_0}\left[e^{f(X_T^{\eps})/\eps}\right] = \sup_{\{x_t\}: x_0 \text{ fixed}} \left[f(x_T) - \frac{1}{2}\int_0^{T}|\dot{x}_{t} - b(t,x_{t})|^2\, dt \right].
\end{equation}
\end{propbox}
This recovers the LDP of Section~\ref{sec:LDP_trajectories}, with Lagrangian $L(t,x, \dot{x}) = \frac{1}{2}|\dot{x} - b(t,x)|^2$. The optimal path---the one achieving the supremum---is the most probable way to realize the rare event, and corresponds to a specific ``optimal noise realization.''

While this path-integral derivation is heuristic, it provides valuable intuition: rare events occur when the noise conspires to push the system along an atypical trajectory, and the probability of such a conspiracy is controlled by the action $\int L\,dt$.

\section{Background on Markov jump processes}
\label{appendix:MJP}

While the main text focuses on diffusion processes, many systems---particularly in chemical kinetics and population dynamics---are better described by Markov jump processes (MJPs). This appendix summarizes the key properties that parallel the SDE theory developed in Appendix~\ref{appendix:ITO}.

\subsection{Definition via transition rates}

An MJP on state space $\Omega \subseteq \mathbb{R}^d$ is specified by $R$ reaction channels with rates $\lambda_j(x) \geq 0$ and displacement vectors $\nu_j \in \mathbb{R}^d$, $j = 1, \ldots, R$. In an infinitesimal time $\Delta t$, the process jumps from state $x$ to state $x + \nu_j$ with probability $\lambda_j(x)\,\Delta t + o(\Delta t)$, and remains at $x$ with probability $1 - \Lambda(x)\,\Delta t + o(\Delta t)$, where $\Lambda(x) = \sum_{j=1}^R \lambda_j(x)$ is the total rate.

This definition has a transparent interpretation in terms of the Gillespie algorithm~\cite{Gillespie1977}, which simulates MJPs exactly. Given the current state $X_t$, the process waits for an exponentially distributed time $\tau \sim \text{Exp}(\Lambda(X_t))$. When the jump occurs, reaction $j$ is selected with probability $\lambda_j(X_t)/\Lambda(X_t)$, and the state updates to $X_{t+\tau} = X_t + \nu_j$. This cycle then repeats from the new state. The memoryless property of the exponential distribution ensures that the process is Markov.

\subsection{The generator}

The generator of an MJP, analogous to~\eqref{eq:generator_def_app} for diffusions, is
\begin{equation}
    (L\phi)(x) = \lim_{\Delta t \downarrow 0} \frac{1}{\Delta t} \mathbb{E}^{X_t=x}[\phi(X_{t+\Delta t}) - \phi(X_t)] = \sum_{j=1}^{R} \lambda_j(x)\left[\phi(x + \nu_j) - \phi(x)\right].
    \label{eq:generator_MJP_app}
\end{equation}
This follows directly from the transition probabilities: in time $\Delta t$, the expected change in $\phi$ is
\begin{equation}
    \mathbb{E}^{X_t=x}[\phi(X_{t+\Delta t}) - \phi(x)] = \sum_{j=1}^R \lambda_j(x)\,\Delta t\, [\phi(x+\nu_j) - \phi(x)] + o(\Delta t).
\end{equation}

\subsection{Backward Kolmogorov equation}

The BKE holds for MJPs with the same interpretation as for SDEs: if $u(t,x)$ satisfies
\begin{equation}
    \partial_t u + Lu = 0, \qquad u(T,x) = F(x),
\end{equation}
then $u(t,x) = \mathbb{E}^{X_t=x}[F(X_T)]$. 

The proof parallels the diffusive case. Define the process $M_t = u(t, X_t)$. Between jumps, $X_t$ is constant, so $dM_t = \partial_t u\,dt$. At a jump time $\tau$ where $X_{\tau^-} = x$ and $X_\tau = x + \nu_j$, the process changes by $u(\tau, x+\nu_j) - u(\tau, x)$. The expected rate of change due to jumps is precisely $(Lu)(t,x)$. Combining these contributions:
\begin{equation}
    \mathbb{E}[dM_t \,|\, X_t = x] = \left(\partial_t u + Lu\right) dt = 0,
\end{equation}
where the last equality uses the BKE. Thus $M_t$ is a martingale, and taking expectations yields $u(0,x_0) = \mathbb{E}[M_0] = \mathbb{E}[M_T] = \mathbb{E}[F(X_T)]$.

\subsection{Forward Kolmogorov equation (Master equation)}

The probability mass function $p(t,x) = \mathbb{P}(X_t = x)$ evolves according to the \textit{master equation}:
\begin{equation}
    \partial_t p(t,x) = \sum_{j=1}^R \left[\lambda_j(x - \nu_j)\,p(t, x-\nu_j) - \lambda_j(x)\,p(t,x)\right].
    \label{eq:master_equation}
\end{equation}
This is the adjoint of the generator, describing probability flow into and out of each state: the first term represents inflow from states $x - \nu_j$ that can jump to $x$, while the second term represents outflow from $x$ to other states.

\section{Variational problems with terminal payoff}
\label{appendix:euler_lagrange}

This appendix derives the Euler--Lagrange and Hamilton equations for variational problems with fixed initial position and a terminal payoff function, yielding the boundary conditions used throughout the main text.

\subsection{Euler--Lagrange equations}

Consider the action
\begin{equation}
    S_T[x] = \int_{0}^{T} L(t,x_t, \dot{x}_t)\, dt - f(x_{T}),
    \label{eq:action-terminal}
\end{equation}
to be minimized over paths with fixed initial condition $x_{t=0} = x_0$ and free final condition $x_{t=T}$. Let $x_t$ be an extremum and $\tilde{x}_t = x_t + \varepsilon\, \delta x_t$ a perturbation with $\delta x_{0} = 0$. The stationarity condition $\frac{d}{d\varepsilon} S_T[\tilde{x}]\big|_{\varepsilon=0} = 0$ gives
\begin{equation}
    0 = \int_{0}^{T}  \left( \frac{\partial L}{\partial x_t}\, \delta x_t + \frac{\partial L}{\partial \dot{x}_t}\, \delta \dot{x}_t \right) dt - \nabla f(x_{T}) \cdot \delta x_{T}.
\end{equation}
Integrating by parts and using $\delta x_{0} = 0$:
\begin{equation}
    0 = \int_{0}^{T}  \left( \frac{\partial L}{\partial x_t} - \frac{d}{dt} \frac{\partial L}{\partial \dot{x_t}} \right) \delta x_t dt + \left( \frac{\partial L}{\partial \dot{x}_t}\bigg|_{T} - \nabla f(x_{T}) \right) \cdot \delta x_{T}.
\end{equation}
Since $\delta x$ is arbitrary, we obtain the Euler--Lagrange equation
\begin{equation}
    \frac{d}{dt} \frac{\partial L}{\partial \dot{x}_t} = \frac{\partial L}{\partial x_t},
    \label{eq:EL}
\end{equation}
with boundary conditions
\begin{align}
    x_{t=0} &= x_0, \label{eq:BC-initial} \\
    \frac{\partial L}{\partial \dot{x}}\bigg|_{T} &= \nabla f(x_{t=T}). \label{eq:BC-final}
\end{align}

\subsection{Hamilton's equations}

The Hamiltonian formulation follows from the Legendre transform. Define the conjugate momentum
\begin{equation}
    \theta_t = \frac{\partial L}{\partial \dot{x}_t},
    \label{eq:momentum-def}
\end{equation}
and the Hamiltonian
\begin{equation}
    H(t,x_t, \theta_t) = \theta_t \cdot \dot{x}_t - L(t,x_t, \dot{x}_t),
    \label{eq:Hamiltonian-def}
\end{equation}
where $\dot{x}_t(t,x_t, \theta_t)$ is obtained by inverting the momentum relation~\eqref{eq:momentum-def}. Taking differentials of~\eqref{eq:Hamiltonian-def}:
\begin{equation}
    dH = \dot{x}_t \cdot d\theta_t + \theta_t \cdot d\dot{x}_t - \frac{\partial L}{\partial t}\, dt - \frac{\partial L}{\partial x_t} \cdot dx_t - \frac{\partial L}{\partial \dot{x}_t} \cdot d\dot{x}_t.
\end{equation}
Using $\theta_t = \partial L / \partial \dot{x}_t$, the $d\dot{x}_t$ terms cancel, yielding
\begin{equation}
    dH = \dot{x}_t \cdot d\theta_t - \frac{\partial L}{\partial t}\, dt - \frac{\partial L}{\partial x_t} \cdot dx_t 
\end{equation}
Comparing with $dH = \frac{\partial H}{\partial x_t} \cdot dx_t + \frac{\partial H}{\partial \theta_t} \cdot d\theta_t + \frac{\partial H}{\partial t}\, dt$, we identify
\begin{equation}
    \frac{\partial H}{\partial \theta_t} = \dot{x}_t, \qquad \frac{\partial H}{\partial x_t} = -\frac{\partial L}{\partial x_t}, \qquad \frac{\partial H}{\partial t} = -\frac{\partial L}{\partial t}.
\end{equation}
The Euler--Lagrange equation~\eqref{eq:EL} becomes $\dot{\theta}_t = \partial L / \partial x_t = -\partial H / \partial x_t$, yielding Hamilton's equations:
\begin{equation}
    \dot{x}_t = \frac{\partial H}{\partial \theta_t}, \qquad \dot{\theta}_t = -\frac{\partial H}{\partial x_t}.
    \label{eq:Hamilton}
\end{equation}

The boundary conditions~\eqref{eq:BC-initial}--\eqref{eq:BC-final} become:
\begin{equation}
    x_{t=0} = x_0, \qquad \theta_{t=T} = \nabla f(x_{T}).
    \label{eq:BC-Hamiltonian}
\end{equation}
These are the boundary condition~\eqref{eq:boundary_conditions} in the main text: position is fixed at the initial time, while momentum is determined by the terminal payoff at the final time.

\begin{remark}[Action in Hamiltonian form]
Along an extremal path, the action~\eqref{eq:action-terminal} can be rewritten using~\eqref{eq:Hamiltonian-def}:
\begin{equation}
    S_T[x^*] = \int_{0}^{T} \left[ \theta \cdot \dot{x}_t - H(x_t, \theta_t, t) \right] dt - f(x_{T}).
\end{equation}
This is the form used for minimum action methods, where the saddle-point structure (minimization over $x$, maximization over $\theta$) is exploited numerically.
\end{remark}
\bibliography{ref_eric}

@article{yan_learning_2022,
	title = {Learning nonequilibrium control forces to characterize dynamical phase transitions},
	volume = {105},
	url = {https://link.aps.org/doi/10.1103/PhysRevE.105.024115},
	doi = {10.1103/PhysRevE.105.024115},
	number = {2},
	urldate = {2023-09-05},
	journal = {Physical Review E},
	publisher = {American Physical Society},
	author = {Yan, Jiawei and Touchette, Hugo and Rotskoff, Grant M.},
	month = feb,
	year = {2022},
	pages = {024115},
}

@article{yan_physicsinformed_2022,
	title = {Physics-informed graph neural networks enhance scalability of variational nonequilibrium optimal control},
	volume = {157},
	issn = {0021-9606},
	url = {https://aip-scitation-org.stanford.idm.oclc.org/doi/10.1063/5.0095593},
	doi = {10.1063/5.0095593},
	number = {7},
	urldate = {2022-10-13},
	journal = {The Journal of Chemical Physics},
	publisher = {American Institute of Physics},
	author = {Yan, Jiawei and Rotskoff, Grant M.},
	month = aug,
	year = {2022},
	pages = {074101},
}

@book{olver1974,
  title = {Asymptotics and {{Special Functions}}},
  author = {Olver, F.W.J.},
  year = 1974,
  publisher = {Elsevier},
  doi = {10.1016/B978-0-12-525850-0.50011-0},
  urldate = {2026-02-20},
  copyright = {https://www.elsevier.com/tdm/userlicense/1.0/},
  isbn = {978-0-12-525850-0}
}

@book{dembo2009large,
  title={Large Deviations Techniques and Applications},
  author={Dembo, Amir and Zeitouni, Ofer},
  edition={2nd},
  year={2009},
  publisher={Springer},
  series={Stochastic Modelling and Applied Probability},
  volume={38}
}

@book{varadhan1984large,
  title={Large Deviations and Applications},
  author={Varadhan, S. R. Srinivasa},
  year={1984},
  publisher={SIAM}
}

@book{FreidlinWentzell2012,
  author    = {Freidlin, Mark I. and Wentzell, Alexander D.},
  title     = {Random Perturbations of Dynamical Systems},
  edition   = {3rd},
  publisher = {Springer},
  address   = {Berlin},
  year      = {2012},
  series    = {Grundlehren der mathematischen Wissenschaften},
  volume    = {260}
}

@book{Oksendal2003,
  author    = {{\O}ksendal, Bernt},
  title     = {Stochastic Differential Equations: An Introduction with Applications},
  edition   = {6th},
  publisher = {Springer},
  address   = {Berlin},
  year      = {2003}
}

@article{bouchet2016large,
  title={Large deviations in fast--slow systems},
  author={Bouchet, Freddy and Grafke, Tobias and Tangarife, Tom{\'a}s and Vanden-Eijnden, Eric},
  journal={Journal of Statistical Physics},
  volume={162},
  number={4},
  pages={793--812},
  year={2016},
  publisher={Springer}
}

@article{Cloning1,
  doi = {10.1088/1742-5468/2007/03/P03004},
  year = {2007},
  volume = {2007},
  number = {03},
  pages = {P03004},
  author = {Lecomte, Vivien and Tailleur, Julien},
  title = {A numerical approach to large deviations in continuous time},
  journal = {Journal of Statistical Mechanics: Theory and Experiment}
}

@article{Cloning2,
  author={Giardina, Cristian and Kurchan, Jorge and Lecomte, Vivien and Tailleur, Julien},
  title={Simulating Rare Events in Dynamical Processes},
  journal={Journal of Statistical Physics},
  year={2011},
  volume={145},
  number={4},
  pages={787-811},
  doi={10.1007/s10955-011-0350-4}
}

@article{schlogl1972chemical,
  title={Chemical reaction models for non-equilibrium phase transitions},
  author={Schl{\"o}gl, Friedrich},
  journal={Zeitschrift f{\"u}r Physik},
  volume={253},
  number={2},
  pages={147--161},
  year={1972},
  publisher={Springer}
}

@article{PhysRevE.86.040103,
  title = {Exchange of stability as a function of system size in a nonequilibrium system},
  author = {T\u{a}nase-Nicola, Sorin and Lubensky, David K.},
  journal = {Physical Review E},
  volume = {86},
  issue = {4},
  pages = {040103},
  year = {2012},
  doi = {10.1103/PhysRevE.86.040103}
}

@article{Gillespie1977,
  author  = {Gillespie, Daniel T.},
  title   = {Exact stochastic simulation of coupled chemical reactions},
  journal = {Journal of Physical Chemistry},
  volume  = {81},
  pages   = {2340--2361},
  year    = {1977},
  doi     = {10.1021/j100540a008}
}

@article{PhysRevB.66.052301,
  title = {String method for the study of rare events},
  author = {E, Weinan and Ren, Weiqing and Vanden-Eijnden, Eric},
  journal = {Physical Review B},
  volume = {66},
  issue = {5},
  pages = {052301},
  year = {2002},
  doi = {10.1103/PhysRevB.66.052301}
}

@article{ren2007simplified,
  title={Simplified and improved string method for computing the minimum energy paths in barrier-crossing events},
  author={Ren, Weiqing and Vanden-Eijnden, Eric},
  journal={The Journal of Chemical Physics},
  volume={126},
  number={16},
  pages={164103},
  year={2007},
  publisher={AIP Publishing}
}

@article{vanden2008geometric,
  title={The geometric minimum action method for computing minimum energy paths},
  author={Vanden-Eijnden, Eric and Heymann, Matthias},
  journal={The Journal of Chemical Physics},
  volume={128},
  number={6},
  pages={061103},
  year={2008},
  publisher={AIP Publishing}
}

@article{zakine2023minimum,
  title = {Minimum-Action Method for Nonequilibrium Phase Transitions},
  author = {Zakine, Ruben and Vanden-Eijnden, Eric},
  journal = {Physical Review X},
  volume = {13},
  issue = {4},
  pages = {041044},
  year = {2023},
  doi = {10.1103/PhysRevX.13.041044}
}

@article{grafke2019numerical,
  title={Numerical computation of rare events via large deviation theory},
  author={Grafke, Tobias and Vanden-Eijnden, Eric},
  journal={Chaos: An Interdisciplinary Journal of Nonlinear Science},
  volume={29},
  number={6},
  pages={063118},
  year={2019},
  publisher={AIP Publishing}
}

@article{nn_am, 
  author = {Boffi, Nicholas M. and Vanden-Eijnden, Eric},
  title = {Deep learning probability flows and entropy production rates in active matter},
  journal = {Proceedings of the National Academy of Sciences},
  volume = {121},
  number = {25},
  pages = {e2318106121},
  year = {2024},
  doi = {10.1073/pnas.2318106121}
}

@article{boffi2023probability,
  title={Probability flow solution of the Fokker--Planck equation},
  author={Boffi, Nicholas M and Vanden-Eijnden, Eric},
  journal={Machine Learning: Science and Technology},
  volume={4},
  number={3},
  pages={035012},
  year={2023},
  publisher={IOP Publishing}
}

@article{simonnet2023computing,
  title={Computing non-equilibrium trajectories by a deep learning approach},
  author={Simonnet, Eric},
  journal={Journal of Computational Physics},
  volume={491},
  pages={112349},
  year={2023},
  publisher={Elsevier}
}

@book{ELiVE2019,
  author    = {E, Weinan and Li, Tiejun and Vanden-Eijnden, Eric},
  title     = {Applied Stochastic Analysis},
  publisher = {American Mathematical Society},
  address   = {Providence, RI},
  year      = {2019},
  series    = {Graduate Studies in Mathematics},
  volume    = {199},
  isbn      = {978-1-4704-4933-9}
}

@book{Gardiner2009,
  author    = {Gardiner, Crispin W.},
  title     = {Stochastic Methods: A Handbook for the Natural and Social Sciences},
  edition   = {4th},
  publisher = {Springer},
  address   = {Berlin},
  year      = {2009},
  series    = {Springer Series in Synergetics},
  isbn      = {978-3-540-70712-7}
}

@article{schuttler2024,
  author = {Sch\"uttler, Janik and Jack, Robert L. and Cates, Michael E.},
  title = {Nonreciprocal dynamics with weak noise: aperiodic ``{E}scher cycles'' and their quasipotential landscape},
  journal = {arXiv preprint},
  year = {2024},
  eprint = {2512.08210},
  archivePrefix = {arXiv},
  primaryClass = {cond-mat.stat-mech}
}

@article{vanden2012rare,
  title={Rare event simulation of small noise diffusions},
  author={Vanden-Eijnden, Eric and Weare, Jonathan},
  journal={Communications on Pure and Applied Mathematics},
  volume={65},
  number={12},
  pages={1770--1803},
  year={2012},
  publisher={Wiley Online Library}
}

@article{varadhan1966,
  author = {Varadhan, S. R. S.},
  title = {Asymptotic probabilities and differential equations},
  journal = {Communications on Pure and Applied Mathematics},
  volume = {19},
  number = {3},
  pages = {261--286},
  year = {1966},
  publisher = {Wiley}
}

@article{gartner1977,
  author = {G{\"a}rtner, J{\"u}rgen},
  title = {On large deviations from the invariant measure},
  journal = {Theory of Probability \& Its Applications},
  volume = {22},
  number = {1},
  pages = {24--39},
  year = {1977},
  publisher = {SIAM}
}

@article{ellis1984,
  author = {Ellis, Richard S.},
  title = {Large deviations for a general class of random vectors},
  journal = {The Annals of Probability},
  volume = {12},
  number = {1},
  pages = {1--12},
  year = {1984},
  publisher = {Institute of Mathematical Statistics}
}

@article{schilder1966,
  author = {Schilder, M.},
  title = {Some asymptotic formulas for {W}iener integrals},
  journal = {Transactions of the American Mathematical Society},
  volume = {125},
  number = {1},
  pages = {63--85},
  year = {1966},
  publisher = {American Mathematical Society}
}

@article{touchette2009,
  author = {Touchette, Hugo},
  title = {The large deviation approach to statistical mechanics},
  journal = {Physics Reports},
  volume = {478},
  number = {1-3},
  pages = {1--69},
  year = {2009},
  publisher = {Elsevier}
}

@article{cramer1938,
  author = {Cram{\'e}r, Harald},
  title = {Sur un nouveau th{\'e}or{\`e}me-limite de la th{\'e}orie des probabilit{\'e}s},
  journal = {Actualit{\'e}s Scientifiques et Industrielles},
  volume = {736},
  pages = {5--23},
  year = {1938},
  publisher = {Hermann, Paris}
}

@book{denhollander2000,
  author = {den Hollander, Frank},
  title = {Large Deviations},
  series = {Fields Institute Monographs},
  volume = {14},
  publisher = {American Mathematical Society},
  address = {Providence, RI},
  year = {2000}
}

@book{fengkurtz2006,
  author = {Feng, Jin and Kurtz, Thomas G.},
  title = {Large Deviations for Stochastic Processes},
  series = {Mathematical Surveys and Monographs},
  volume = {131},
  publisher = {American Mathematical Society},
  address = {Providence, RI},
  year = {2006}
}

@book{ellis1985,
  author = {Ellis, Richard S.},
  title = {Entropy, Large Deviations, and Statistical Mechanics},
  series = {Grundlehren der mathematischen Wissenschaften},
  volume = {271},
  publisher = {Springer},
  address = {New York},
  year = {1985}
}

@book{flemingsoner2006,
  author = {Fleming, Wendell H. and Soner, H. Mete},
  title = {Controlled {M}arkov Processes and Viscosity Solutions},
  edition = {2nd},
  series = {Stochastic Modelling and Applied Probability},
  volume = {25},
  publisher = {Springer},
  address = {New York},
  year = {2006}
}

@article{crandallishiilions1992,
  author = {Crandall, Michael G. and Ishii, Hitoshi and Lions, Pierre-Louis},
  title = {User's guide to viscosity solutions of second order partial differential equations},
  journal = {Bulletin of the American Mathematical Society},
  volume = {27},
  number = {1},
  pages = {1--67},
  year = {1992}
}

@article{kurtz1972,
  author = {Kurtz, Thomas G.},
  title = {The relationship between stochastic and deterministic models for chemical reactions},
  journal = {The Journal of Chemical Physics},
  volume = {57},
  number = {7},
  pages = {2976--2978},
  year = {1972}
}

@book{ethierkurtz1986,
  author = {Ethier, Stewart N. and Kurtz, Thomas G.},
  title = {Markov Processes: Characterization and Convergence},
  publisher = {Wiley},
  address = {New York},
  year = {1986}
}

@book{kipnislandim1999,
  author = {Kipnis, Claude and Landim, Claudio},
  title = {Scaling Limits of Interacting Particle Systems},
  series = {Grundlehren der mathematischen Wissenschaften},
  volume = {320},
  publisher = {Springer},
  address = {Berlin},
  year = {1999}
}

@article{donskervaradhan1975,
  author = {Donsker, Monroe D. and Varadhan, S. R. S.},
  title = {Asymptotic evaluation of certain {M}arkov process expectations for large time, {I}},
  journal = {Communications on Pure and Applied Mathematics},
  volume = {28},
  number = {1},
  pages = {1--47},
  year = {1975}
}

@article{chetritetouchette2015,
  author = {Chetrite, Rapha{\"e}l and Touchette, Hugo},
  title = {Nonequilibrium {M}arkov processes conditioned on large deviations},
  journal = {Annales Henri Poincar{\'e}},
  volume = {16},
  number = {9},
  pages = {2005--2057},
  year = {2015}
}

@article{weinan2004mam,
  author = {E, Weinan and Ren, Weiqing and Vanden-Eijnden, Eric},
  title = {Minimum action method for the study of rare events},
  journal = {Communications on Pure and Applied Mathematics},
  volume = {57},
  number = {5},
  pages = {637--656},
  year = {2004}
}

@article{bouchet2016perturbative,
  author = {Bouchet, Freddy and Gaw{\c{e}}dzki, Krzysztof and Nardini, Cesare},
  title = {Perturbative calculation of quasi-potential in non-equilibrium diffusions: a mean-field example},
  journal = {Journal of Statistical Physics},
  volume = {163},
  number = {5},
  pages = {1157--1210},
  year = {2016}
}

@article{cameron2012finding,
  author = {Cameron, Maria K.},
  title = {Finding the quasipotential for nongradient {SDE}s},
  journal = {Physica D: Nonlinear Phenomena},
  volume = {241},
  number = {18},
  pages = {1532--1550},
  year = {2012}
}

@article{raissi2019physics,
  author = {Raissi, Maziar and Perdikaris, Paris and Karniadakis, George E.},
  title = {Physics-informed neural networks: A deep learning framework for solving forward and inverse problems involving nonlinear partial differential equations},
  journal = {Journal of Computational Physics},
  volume = {378},
  pages = {686--707},
  year = {2019},
  publisher = {Elsevier}
}

@book{DupuisEllis1997,
    author    = {Dupuis, Paul and Ellis, Richard S.},
    title     = {A Weak Convergence Approach to the Theory of Large Deviations},
    series    = {Wiley Series in Probability and Statistics},
    publisher = {John Wiley \& Sons},
    address   = {New York},
    year      = {1997},
    doi       = {10.1002/9781118165904}
}

@article{NueskenRichter2021,
    author  = {N\"usken, Nikolas and Richter, Lorenz},
    title   = {Solving High-Dimensional {H}amilton--{J}acobi--{B}ellman {PDE}s 
               Using Neural Networks: Perspectives from the Theory of Controlled 
               Diffusions and Measures on Path Space},
    journal = {Partial Differential Equations and Applications},
    volume  = {2},
    number  = {4},
    pages   = {48},
    year    = {2021},
    doi     = {10.1007/s42985-021-00102-x}
}

\end{document}